\documentclass[Afour,sageh,times]{sagej}

\usepackage{moreverb,url}

\usepackage[colorlinks,bookmarksopen,bookmarksnumbered,citecolor=red,urlcolor=red]{hyperref}

\newcommand\BibTeX{{\rmfamily B\kern-.05em \textsc{i\kern-.025em b}\kern-.08em
		T\kern-.1667em\lower.7ex\hbox{E}\kern-.125emX}}

\def\volumeyear{2019}

\usepackage{todonotes}
\usepackage{xspace}
\usepackage{url}
\usepackage{hyperref}
\usepackage{amsmath}
\usepackage{caption} 
\usepackage{subcaption}
\usepackage{amsopn}

\usepackage{lipsum}
\usepackage{amsfonts}
\usepackage{graphicx}
\usepackage{epstopdf}
\usepackage{nameref}

\usepackage{lineno}
\usepackage[algo2e,ruled,vlined,linesnumbered]{algorithm2e}

\newcommand{\subfigwidthw}{70mm}
\newcommand{\subfigwidthww}{130mm}

\newcommand{\MHz}[0]{\,MHz\xspace}

\newcommand{\GB}[0]{\,GB\xspace}

\newcommand{\ylrev}[1]{{\color{black}{#1}}}

\DeclareMathOperator{\diag}{diag}
\DeclareMathOperator*{\argmax}{arg\,max}
\DeclareGraphicsExtensions{.pdf,.jpeg,.png,.eps,.jpg}

\makeatletter
\edef\orig@output{\the\output}
\output{\setbox\@cclv\vbox{\unvbox\@cclv\vspace{0pt plus 10pt}}\orig@output}
\makeatother

\begin{document}
	
	\runninghead{Liu, Sid-Lakhdar, Rebrova, Ghysels and Li}
	
	\title{A Parallel Hierarchical Blocked Adaptive Cross Approximation Algorithm}
	
	\author{Yang Liu\affilnum{1}, Wissam Sid-Lakhdar\affilnum{1}, Elizaveta Rebrova\affilnum{2}, Pieter Ghysels\affilnum{1} and Xiaoye Sherry Li\affilnum{1}}
	
	\affiliation{\affilnum{1}Computational Research Division, Lawrence Berkeley National Laboratory, Berkeley, CA, USA\\
		\affilnum{2}Department of Mathematics, University of California, Los Angeles, CA, USA}
	
	\corrauth{Yang Liu, Computational Research Divisio
		Lawrence Berkeley National Laboratory,
		Berkeley,
		CA,
		USA.}
	
	\email{liuyangzhuan@lbl.gov}
	
	\begin{abstract}
		This paper presents a low-rank decomposition algorithm assuming any matrix element can be computed in $O(1)$ time. The proposed algorithm first computes rank-revealing decompositions of sub-matrices with a blocked adaptive cross approximation (BACA) algorithm, and then applies a hierarchical merge operation via truncated singular value decompositions (H-BACA). The proposed algorithm significantly improves the convergence of the baseline ACA algorithm and achieves reduced computational complexity compared to the full decompositions such as rank-revealing QR. Numerical results demonstrate the efficiency, accuracy and parallel scalability of the proposed algorithm.
	\end{abstract}
	
	
	\maketitle

	\section{Introduction}\label{intro}
	Rank-revealing decomposition algorithms are important numerical linear algebra tools for compressing high-dimensional data, accelerating solution of integral and partial differential equations, constructing efficient machine learning algorithms, and analyzing numerical algorithms, etc, as matrices arising from many science and engineering applications oftentimes exhibit numerical rank-deficiency. Despite the favorable $O(nr)$ memory footprint of such decompositions with $n$ and $r$ respectively denoting the matrix dimension (assuming a square matrix) and \ylrev{the} numerical rank, the computational cost can be expensive. Existing rank-revealing decompositions such as truncated singular value decomposition (SVD), column-pivoted QR (QRCP), CUR decomposition, interpolative decomposition (ID), and rank-revealing LU typically require at least $O(n^2r)$ operations \cite{MingGu1996RRQR,Martinsson2005LR,Martinsson2017IDCUR,Mahoney2014CUR}. 
	This complexity can be reduced to $O(n^2\mathrm{log}\: r + nr^2)$ by structured random matrix projection-based algorithms \cite{Martinsson2017IDCUR, Liberty2007Randomized}. In addition, faster algorithms are available in the following three scenarios. 1. When each element entry can be computed in $O(1)$ CPU time with prior knowledge (i.e., smoothness, sparsity, or leverage scores) about the matrix, faster algorithms such as randomized CUR and adaptive cross approximation (ACA) \cite{Bebendorf2000ACA,Bebendorf2006BEMACA,Zhao2005ACAEM} algorithms can achieve $O(nr^2)$ complexity. However, the robustness of these algorithms relies heavily on matrix properties that are not always present in practice. 
	2. When the matrix can be rapidly applied to arbitrary vectors, algorithms such as randomized SVD, QR and UTV (T lower or upper triangular) \cite{Liberty2007Randomized,MingGu2017qr,MingGu2018SRQR,Martinsson2017utv} can be utilized to achieve quasi-linear complexity. 3. Finally, given a matrix with missing entries, the \ylrev{low-rank} decomposition can be constructed via matrix completion algorithms \cite{Candes2009MatrixCompletion,Laura2010Grouse} in quasi-linear time assuming incoherence properties of the matrices (i.e., projection of natural basis vectors onto the space spanned by singular vectors of the matrix should not be very sparse). This work concerns the development of a practical algorithm, in application scenario 1, that improves the robustness of ACA algorithms while maintaining reduced complexity for broad classes of matrices. 
	
	The partially-pivoted ACA algorithm, closely related to LU with rook pivoting \cite{FOSTER1997ROOK}, constructs an LU-type decomposition upon accessing one row and column per iteration. For matrices resulting from asymptotically smooth kernels, ACA is a rank-revealing and optimal-complexity algorithm that converges in $O(k)$ iterations \cite{Bebendorf2000ACA}. Despite its favorable computational complexity, it is well-known that the ACA algorithm suffers from deteriorated convergence and/or \ylrev{premature} termination for non-smooth, sparse \ylrev{and/or coherent} matrices \cite{Heldring2014ACAconvergence}. Hybrid methods or improved convergence criteria (\ylrev{e.g.}, hybrid ACA-CUR, averaging, statistical norm estimation) have been proposed to partially alleviate the problem \ylrev{\cite{Heldring2015stochasticACA,Grasedyck2005HACA}}. \ylrev{The main difficulty of leveraging ACA as robust algebraic tools for general low-rank matrices results from ACA's partial pivot-search strategy to attain low complexity. In addition to the abovementioned remedies, another possibility to improve ACA's robustness is to search for pivots in a wider range of rows/columns without sacrificing too much computational efficiency. Here we consider two different strategies: 1. Instead of searching one row/column per iteration as in ACA, it is possible to search a block of rows/columns to find multiple pivots together. 2. Instead of applying ACA directly on the entire matrix, it is possible to start with compressing submatrices via ACA and then merge the results as one low-rank product. In extreme cases (e.g., when block size equals matrix dimension or submatrix dimension equals one), these strategies lead to quadratic computational costs. Therefore, it is valuable to address the question: for what matrix kernels and under what block/submatrix sizes will these strategies retain low complexity.}   
	
	\ylrev{For the first strategy}, this work proposes a blocked ACA algorithm (BACA) that extracts a block row/column per iteration to significantly improve convergence of the baseline ACA algorithms. The blocked version also enjoys higher flop performance as it involves mainly BLAS-3 operations. Compared to the aforementioned remedies, the proposed algorithm provides a unified framework to balance robustness and efficiency. Upon increasing the block size \ylrev{(i.e., the number of rows/columns per iteration)}, the algorithm gradually changes from ACA to ID. \ylrev{For the second strategy, the proposed algorithm further subdivides the matrix} into $n_b$ submatrices compressed via BACA, followed by a hierarchical merge algorithm \ylrev{leveraging low-rank arithmetic} \ylrev{\cite{Hackbusch2002Hmatrix,Hackbusch2003Hmatrix}}. The overall cost of this H-BACA algorithm is at most $O(\sqrt{n_b}nr^2)$ \ylrev{assuming the block size in BACA is less than the rank}. In other words, the proposed H-BACA algorithm is a general numerical linear algebra tool as an alternative to ACA, SVD, QR, etc.
	In addition, the overall algorithm can be parallelized using distributed-memory linear algebra packages such as ScaLAPACK \cite{scalapack} which avoids the difficulty of efficient parallelization of plain ACA algorithms. Numerical results illustrate good accuracy, efficiency and parallel performance. In addition, the proposed algorithm can be used as a general low-rank compression tool for constructing hierarchical matrices \cite{Liza2018Clustering}.   
	
	\section{Notation}\label{notation}
	Throughout this paper, we adopt the Matlab notation of matrices and vectors.
	Submatrices of a matrix $A$ are denoted $A(I,J)$, $A(:,J)$ or $A(I,:)$ where $I$, $J$ are index sets. Similarly, subvectors of a column vector $u$ are denoted $u(I)$. An index set $I$ permuted by $J$ reads $I(J)$. Transpose, inverse, pseudo-inverse of $A$ are $A^t$, $A^{-1}$, $A^\dagger$. $\left\lVert A \right\rVert_F$ and $\left\lVert u \right\rVert_2$ denote Frobenius norm and 2-norm. \ylrev{Note that $u$ refers to a $n\times 1$ column vector}. Vertical and horizontal concatenations of $A$, $B$ are $[A;B]$ and $[A,B]$. \ylrev{Element-wise multiplication of $A$ and $B$ is $A\circ B$}. All matrices are real-valued unless otherwise stated. It is assumed for $A\in\mathbb{R}^{m\times n}$, $m=O(n)$, but the proposed algorithms also apply to complex-valued and tall-skinny / short-fat matrices. We denote truncated SVD as $[U,\Sigma,V,r]=\mathtt{SVD}(A,\epsilon)$ with $U\in\mathbb{R}^{m\times r}$, $V^t\in\mathbb{R}^{n\times r}$ column orthogonal, $\Sigma\in\mathbb{R}^{r\times r}$ diagonal, and $r$ being $\epsilon$-rank \ylrev{defined by $r=\min\{k\in\mathbb{N}: \Sigma_{k+1,k+1}<\epsilon\Sigma_{1,1}\}$}. We denote QRCP as $[Q,T,J]=\mathtt{QR}(A,r)$ or $[Q,T,J]=\mathtt{QR}(A,\epsilon)$ with $Q\in\mathbb{R}^{m\times r}$ column orthogonal, $T\in\mathbb{R}^{r\times n}$ upper triangular, $J$ being column pivots,
	and $\epsilon$ and $r$ being the prescribed accuracy and rank, respectively.
	QR without column-pivoting is simply written as $[Q,T]=\mathtt{QR}(A)$. \ylrev{Cholesky decomposition without pivoting is written as $T=\mathtt{Chol}(A)$ with $T$ upper triangular. $\mathrm{log}n$ means logarithm of $n$ to the base 2.}
	
	\section{Algorithm Description}\label{algo_des}
	\subsection{Adaptive Cross Approximation}
	Before describing the proposed algorithm, we first briefly summarize the baseline ACA algorithm \cite{Zhao2005ACAEM}. Consider a matrix $A\in\mathbb{R}^{m\times n}$ of $\epsilon$-rank $r$, 
	the ACA algorithm approximates $A$ by a sequence of rank-1 outer-products as 
	\begin{linenomath*}
		\begin{align}
		A \approx UV=\sum_{k=1}^{r}u_kv_k^t 
		\end{align}  
	\end{linenomath*}
	At each iteration $k$, the algorithm selects column $u_k$ (pivot $j_k$ \ylrev{from remaining columns}) and row $v_k^t$ (pivot $i_k$ \ylrev{from remaining rows}) from the residual matrix $E_{k-1}=A-\sum_{i=1}^{k-1}u_iv_i^t$ corresponding to an element \ylrev{denoted by} $E_{k-1}(i_k,j_k)$ with sufficiently large magnitude. \ylrev{Note that $u_k$ and $v_k$ are $m\times 1$ and $n\times 1$ vectors.} The partially-pivoted ACA algorithm (ACA for short), selecting \ylrev{$j_k, i_k$} by only looking at previously selected rows and columns, is described as Algorithm \ref{algo_aca}. Specifically, each iteration $k$ selects pivot $i_k$ used in the current iteration and pivot $j_{k+1}$ for the next iteration (via line \ref{ik} and \ref{jk}) as
	\begin{linenomath*}
		\begin{align}
		\ylrev{i_k=\argmax_{i\neq i_1,...,i_{k-1}} \lvert E_{k-1}(:,j_k)\rvert}\label{iik}\\
		\ylrev{j_{k+1}=\argmax_{j\neq j_1,...,j_{k}} \lvert E_{k-1}(i_{k},:)\rvert}\label{jjk}
		\end{align}  
	\end{linenomath*}
	\ylrev{and $j_1$ is a random initial column index.} Note that $i_k\neq i_1,...,i_{k-1}$ and $j_k\neq j_1,...,j_{k-1}$ are enforced. The iteration is terminated when $\nu<\epsilon\mu$ with
	\begin{linenomath*}
		\begin{align}
		\nu=\left\lVert u_kv_k^t \right\rVert_F\approx\left\lVert A-UV \right\rVert_F,~~\mu=\left\lVert UV \right\rVert_F\approx\left\lVert A \right\rVert_F \label{aca_stop}
		\end{align}  
	\end{linenomath*}
	and $\epsilon$ is the prescribed tolerance. Note that each iteration requires only $O(nr_k)$ flop operations with $r_k$ denoting currently revealed numerical rank. The overall complexity of partially-pivoted ACA scales as $O(nr^2)$ when the algorithm converges in $O(r)$ iterations. Despite the favorable complexity, the convergence of ACA \ylrev{for general rank-deficient matrices} is unsatisfactory. For many rank-deficient matrices arising from the numerical solution of PDEs, signal processing and data science, ACA oftentimes either requires $O(n)$ iterations or exhibits \ylrev{premature} termination. First, as ACA does not search the full residual matrices for the largest element, it cannot avoid selection of smaller pivots for general rank-deficient matrices and may require $O(n)$ iterations. Second, the approximation $\left\lVert u_kv_k^t \right\rVert_F$ in (\ref{aca_stop}) often causes the premature termination with the selection of smaller pivots. Remedies such as averaged stopping criteria \cite{Zhou2017upgradedACA}, stochastic error estimation \cite{Heldring2015stochasticACA}, ACA+ \ylrev{\cite{Grasedyck2005HACA}}, and hybrid ACA \ylrev{\cite{Grasedyck2005HACA}} have been developed but they do not generalize to a broad range of applications.

	\IncMargin{1em}
	\begin{algorithm2e}
		\SetKwData{Left}{left}\SetKwData{This}{this}\SetKwData{Up}{up}\SetKwData{Conv}{conv}\SetKwData{Err}{er}
		\SetKwFunction{QR}{QR}\SetKwFunction{SVD}{SVD}\SetKwFunction{LRID}{LRID}\SetKwFunction{LRnorm}{LRnorm}\SetKwFunction{LRnormUp}{LRnormUp}
		\SetKwInOut{Input}{input}\SetKwInOut{Output}{output}
		\SetKwProg{Fn}{Function}{}{}
		\Input{Matrix $A\in\mathbb{R}^{m\times n}$, relative tolerance $\epsilon$}
		\Output{Low-rank approximation of $A\approx UV$ with rank $r$}
		$U=0$, $V=0$, $\mu=0$, $r_0=0$, \ylrev{$j_1$} is a random column index\;
		\For{$k= 1$ \KwTo $\min\{m,n\}$}{
			$\ylrev{u_k}=E_{k-1}(:,j_k)=A(:,\ylrev{j_{k}})-UV(:,\ylrev{j_{k}})$\;
			$i_k=\argmax_i \lvert \ylrev{u_k}(i)\rvert$\; \nllabel{ik}
			$u_k\leftarrow \ylrev{u_k/u_k}(i_k)$\;
			$v_k^t=E_{k-1}(i_k,:)=A(i_k,:)-U(i_k,:)V$\;		
			$\ylrev{j_{k+1}}=\argmax_j \lvert v_k(j)\rvert$\; \nllabel{jk}
			$\nu^2=\left\lVert u_k\right\rVert_2^2\left\lVert v_k\right\rVert_2^2$\;
			$\mu^2\leftarrow\mu^2+\nu^2+2\sum_{j=1}^{k-1}{V(j,:)v_ku_k^tU(:,j) }$\;
			$U\leftarrow [U,u_k], V\leftarrow [V;v_k^t], r_k\ylrev{=} r_{k-1}+1$\;
			Terminate if $\nu<\epsilon\mu$.
		}
		\caption{Adaptive cross approximation algorithm (ACA)}\label{algo_aca}
	\end{algorithm2e}\DecMargin{1em}

	\IncMargin{1em}
	\begin{algorithm2e}
		\SetKwData{Left}{left}\SetKwData{This}{this}\SetKwData{Up}{up}\SetKwData{Conv}{conv}\SetKwData{Err}{er}
		\SetKwFunction{QR}{QR}\SetKwFunction{SVD}{SVD}\SetKwFunction{LRID}{LRID}\SetKwFunction{LRnorm}{LRnorm}\SetKwFunction{LRnormUp}{LRnormUp}
		\SetKwInOut{Input}{input}\SetKwInOut{Output}{output}
		\SetKwProg{Fn}{Function}{}{}
		\Input{Matrix $A\in\mathbb{R}^{m\times n}$, block size $d$, relative tolerance $\epsilon$}
		\Output{Low-rank approximation of $A\approx UV$ with rank $r$}
		$U=0$, $V=0$, $r_0=0$, $\mu=0$, \ylrev{$\bar{J}_1$} is a random index set of cardinality $d$\;
		\For{$k= 1$ \KwTo$~\ylrev{\min\{m,n\}}$}{
			\ylrev{$C_k=E_{k-1}(:,J_k)=A(:,J_{k})-UV(:,J_{k})$}\;
			\ylrev{$[Q_k^c,T_k^c,I_k]=\QR(C_k^t,d)$}, $I_k$ denotes selected skeleton rows\;\nllabel{Ik}
			\ylrev{$R_k=E_{k-1}(I_k,:)=A(I_k,:)-U(I_k,:)V$}\;
			\ylrev{$[Q_{k+1}^r,T_{k+1}^r,J_{k+1}]=\QR(R_k,d)$}, $J_{k+1}$ denotes selected skeleton columns\;\nllabel{Jk}
			\ylrev{$W_k=E_{k-1}(I_k,J_k)=A(I_k,J_k)-U(I_k,:)V(:,J_k)$}\;
			\ylrev{$[U_k,V_k,d_k,\ylrev{\bar{J}}]=\LRID(C_k,W_k,R_k)$}\;\nllabel{ID}
			\ylrev{$I_k\leftarrow I_k([1,d_k]), J_k\leftarrow J_k(\bar{J})$}\;
			\ylrev{$r_k= r_{k-1}+d_k$}\;
			\ylrev{$\nu=\LRnorm(U_k,V_k)$}\;\nllabel{normUVk}
			\ylrev{$\mu\leftarrow\LRnormUp(U,V,\mu,U_k,V_k,\nu)$}\;\nllabel{normUV}
			\ylrev{$U\leftarrow [U,U_k], V\leftarrow [V;V_k]$}\;
			\ylrev{Terminate if $\nu<\epsilon\mu$}.
		}
		\Fn{\LRID $($$C$,$W$,$R$,$\epsilon$$)$}{
			\Input{$C=A(:,J)$, $R=A(I,:)$, $W=A(I,J)$ with $I,J$ of same cardinality}
			\Output{$A\approx UV$ with $U\in\mathbb{R}^{m\times r}, V\in\mathbb{R}^{r\times n}$}
			
			$[Q,T,\bar{J},r]=\QR(W,\epsilon)$\;
			$U=C(:,\bar{J})$\;
			$V=T^{-1}Q^tR$\;
			\Return $U,V,r,\ylrev{\bar{J}}$
		}
		\Fn{\LRnorm $($$U$,$V$$)$}{
			\Input{$A=UV$}
			\Output{$\left\lVert A\right\rVert_F$}
			\ylrev{$T_1=\mathtt{Chol}(U^tU)$}\;
			\ylrev{$T_2=\mathtt{Chol}(VV^t)$}\;
			\Return $\left\lVert T_1T_2^t\right\rVert_F$\;
		}
		\Fn{\LRnormUp $($$U,V,\nu,\bar{U},\bar{V},\bar{\nu}$$)$}{
			\Input{$U\in\mathbb{R}^{m\times r}$, $V\in\mathbb{R}^{r\times n}$, $\bar{U}\in\mathbb{R}^{m\times \bar{r}}$, $\bar{V}\in\mathbb{R}^{\bar{r}\times n}$, $\nu=\left\lVert UV\right\rVert_F$, $\bar{\nu}=\left\lVert \bar{U}\bar{V}\right\rVert_F$}
			\Output{$\left\lVert [U,\bar{U}][V;\bar{V}]\right\rVert_F$}		
			\ylrev{$s=\nu^2+\bar\nu^2+2\sum_{i=1}^{r}\sum_{j=1}^{\bar{r}}{\tilde{V}(i,j)}$ with $\tilde{V}=(V\bar{V}^t)\circ(U^t\bar{U})$}\;
			\Return $\sqrt{s}$	
		}
		\caption{Blocked adaptive cross approximation algorithm (BACA)}\label{algo_baca}
	\end{algorithm2e}\DecMargin{1em}

	\subsection{Blocked Adaptive Cross Approximation}\label{baca_alg}
	Instead of selecting only one column and row from the residual matrix in each ACA iteration, we can select a fixed-size block of columns and rows per iteration to improve the convergence and accuracy of ACA. In addition, many BLAS-1 and BLAS-2 operations of ACA become BLAS-3 operations and hence higher flop performance can be achieved. 
	
	Specifically, the proposed BACA algorithm factorizes $A$  
	\begin{linenomath*}
		\begin{align}
		A \approx UV=\sum_{k=1}^{\ylrev{n_d}}U_kV_k 
		\end{align}  
	\end{linenomath*}
	where $U_k\in\mathbb{R}^{m\times d_k}$ and $V_k\in\mathbb{R}^{d_k\times n}$. In principle, the algorithm selects a block of $d$ rows and columns via cross approximations in the residual matrix and then $d_k\leq d$ ones via rank-revealing algorithms to form a low-rank update at iteration $k$. The total number of iterations is approximately $n_d\approx \lceil r/d \rceil$  if $d_k\approx d$. Instead of selecting row/column pivots via lines \ref{ik} and \ref{jk} of Algorithm \ref{algo_aca}, the proposed algorithm selects row and column index sets $I_k$ and $J_k$ by performing QRCP on $d$ columns (more precisely their transpose) and rows of the residual matrices. This proposed strategy is described in Algorithm \ref{algo_baca}.   
	
\begin{figure}[!tb]
		\centering	
		\begin{subfigure}{\subfigwidthw}
			\includegraphics[width=\subfigwidthw]{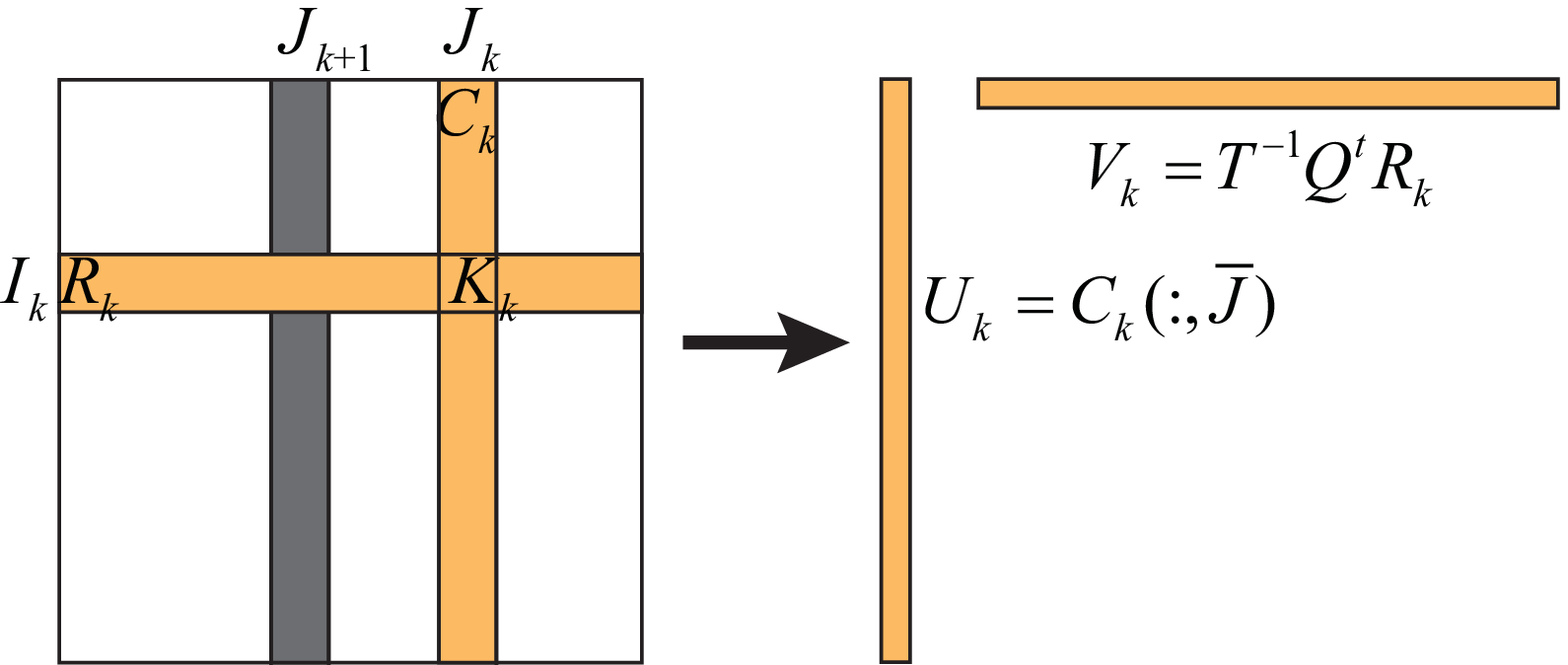}
			\vspace{-20pt}			
			\caption{}\label{fig:IJ}				
		\end{subfigure}
		\begin{subfigure}{\subfigwidthw}
			\includegraphics[width=\subfigwidthw]{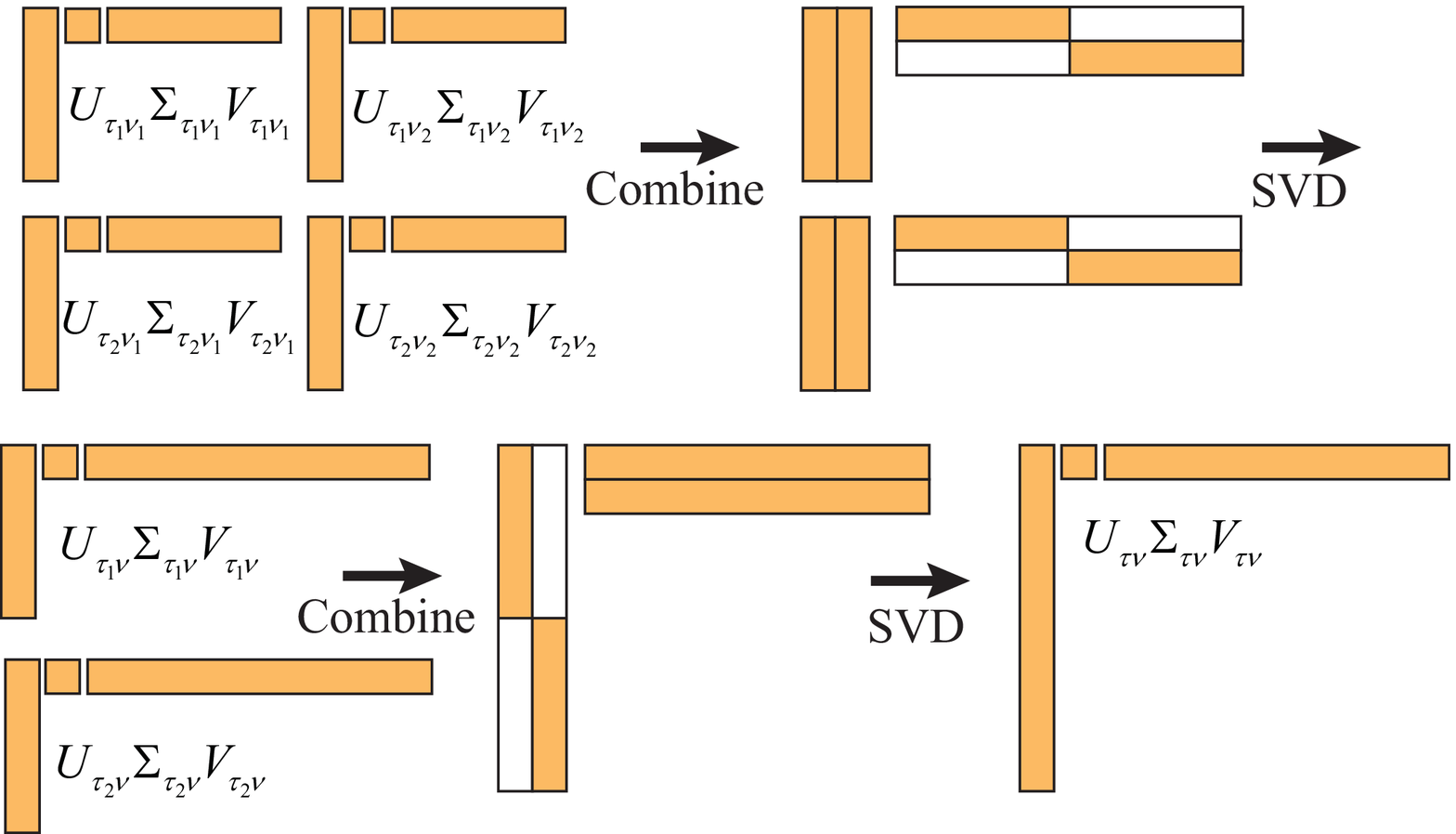}
			\vspace{-20pt}			
			\caption{}\label{fig:merge}				
		\end{subfigure}
		\caption{\ylrev{(a) Selection of $I_k$/$J_k$ and form the low-rank update $U_kV_k$.} (b) Low-rank merge operation} 
	\end{figure}

	Each BACA iteration is composed of three steps. 
	\begin{itemize}
		\item Find block row $I_k$ and block column \ylrev{$J_{k+1}$} by QRCP. \ylrev{Starting with a random column index set $J_{1}$, the block row $I_k$ and the next iteration's block column $J_{k+1}$ are selected by (line \ref{Ik} and \ref{Jk})} 	
		\begin{linenomath*}
			\ylrev{\begin{align} 
				[Q_k^c,T_k^c,I_k]&=\mathtt{QR}(E_{k-1}^t(:,J_k),d) \label{QRIk}\\
				[Q_{k+1}^r,T_{k+1}^r,J_{k+1}]&=\mathtt{QR}(E_{k-1}(I_{k},:),d)\label{QRJk}
				\end{align}}
		\end{linenomath*}
		\ylrev{Here the algorithm first selects $d$ skeleton rows from the submatrix $E_{k-1}(J_k,:)$ (i.e., $d$ columns from its transpose) and then selects $d$ skeleton columns from the submatrix $E_{k-1}(I_{k},:)$ by leveraging the LAPACK implementation of QRCP as it provides a simple way of greedily selecting well-conditioned columns by examining column norms in the $R$ factor at each iteration. Note that many other subset selection algorithms exist in both the machine learning and numerical linear algebra communities (e.g., strong rank-revealing QR \cite{MingGu1996RRQR}, spectrum-revealing QR \cite{2018arXiv180301982F}, and column subset selection problems \cite{Boutsidis2009CSSP}), which ideally pick $d$ matrix columns with maximum volumes. Note that $I_k$ excludes rows selected in previous iterations. To efficiently enforce such condition, the QRCP is performed on the submatrix of $E_{k-1}^t(:,J_k)$ excluding previously selected rows rather than directly on $E_{k-1}^t(:,J_k)$. Similarly, $J_k$ excludes columns selected in previous iterations.} See Fig. \ref{fig:IJ} for an illustration of the procedure. \ylrev{$I_k$ and $J_{k+1}$ are selected by QRCP on the column and transpose of the row marked in yellow, respectively. The column marked in grey is used to select $I_{k+1}$ in the next iteration. For illustration purpose, index sets in Fig. \ref{fig:IJ} consist of contiguous indices.}
		\item Form \ylrev{the factors of} the low-rank product $U_kV_k$. Let $C_k=E_{k-1}(:,J_k)$, $R_k=E_{k-1}(I_k,:)$ and $W_k=E_{k-1}(I_k,J_k)$, $E_{k-1}$ can be approximated by an ID-type decomposition $E_{k-1}\approx C_kW_k^\dagger R_k=U_kV_k$ \cite{Martinsson2017IDCUR} \ylrev{by (\ref{LRID1}) and (\ref{LRID2}). Note that }the pseudo inverse is computed via rank-revealing QR (also see the LRID algorithm at line \ref{ID}). \ylrev{The rank-revealing algorithm is needed as the $d\times d$ block $W_k$ can be further compressed with rank $d_k$. Particularly for matrices where the ACA algorithm tends to fail, the corresponding $d\times d$ matrices $W_k$ in BACA are often rank-deficient. In this case, BACA becomes more robust than ACA as the effective $d_k$ pivots can still be used to generate $d$ columns $J_{k+1}$ for the next iteration (as long as $d_k>0$). Consequently, the effective rank increase is $d_k\leq d$ and the pivot pair $(I_k,J_k)$ is updated in (\ref{IkJkUp}) by the column pivots $\bar{J}$ of QRCP in (\ref{LRID1})}.  
		
		\begin{linenomath*} 
			\begin{align}
			[Q,T,\bar{J}] = \mathtt{QR}(W_k,\epsilon) \mathrm{~with~} Q\in\mathbb{R}^{d\times d_k} \label{LRID1}\\
			U_k = C_k(:,\bar{J}),~V_k = T^{-1}Q^tR_k\label{LRID2}\\
			\ylrev{I_k\leftarrow I_k([1,d_k]), J_k\leftarrow J_k(\bar{J})}\label{IkJkUp}
			\end{align}
		\end{linenomath*}
		\item Compute $\nu=\left\lVert U_kV_k\right\rVert_F$ and update $\mu=\left\lVert UV\right\rVert_F$. Assuming constant block size $d$, the norm of the low-rank update can be computed in $O(nd_k^2)$ operations (line \ref{normUVk}) via
		\begin{linenomath*} 
			\ylrev{\begin{align}
				T_{U_k}&=\mathtt{Chol}(U_k^tU_k), T_{V_k}=\mathtt{Chol}(V_kV_k^t)\label{chol}\\
				\nu&=\left\lVert T_{U_k}T_{V_k}^t\right\rVert_F\label{numu}
				\end{align}} 
		\end{linenomath*}
		Once $\nu$ is computed, the norm of $UV$ can be updated efficiently in $O(nr_kd_k)$ operations (line \ref{normUV}) as
		\begin{linenomath*} 
			\begin{align}
			\ylrev\mu^2\leftarrow\mu^2+\nu^2+2\sum_{i=1}^{r_{k-1}}\sum_{j=1}^{d_k}{\tilde{V}(i,j)}\nonumber\\
			\tilde{V}=(VV_k^t)\circ(U^tU_k)\label{normup}
			\end{align}
		\end{linenomath*} 
		where $r_k$ represents the column dimension of $U$ at iteration $k$. \ylrev{Note that the matrix multiplications in (\ref{chol}) and (\ref{normup}) involving $V_k$ and $V$ (and similarly for those involving $U_k$ and $U$) can be performed as $[V,V_k]V_k^t$ to further improve the computational efficiency.} Then the algorithm updates $U$, $V$ as $[U,U_k]$, $[V;V_k]$ and tests the stopping criteria $\nu<\epsilon\mu$. Note that $\nu,\mu$ with larger $d$ provides better approximations to the exact stop criteria compared to those in (\ref{aca_stop}) hence can significantly reduce the chance of premature termination. 
	\end{itemize}
	
	\ylrev{We would like to highlight the difference between the proposed BACA algorithm and existing ACA algorithms. First, as BACA selects a block of rows and columns per iteration as opposed to a single row and column in the baseline ACA algorithm, the convergence behavior and flop performance can be significantly improved. In the existing ACA algorithms, convergence can also be improved by leveraging averaged stopping criteria \cite{Zhou2017upgradedACA} or searching a single pivot in a broader range of rows and columns (e.g., fully-pivoted ACA). However, they still find one row or column at a time in each iteration and hence suffer from poor flop performance. Moreover, they cannot utilize strong rank revealing algorithms to select skeleton rows and columns with better volume (determinant in modulus) qualities. Second, BACA also has important connections to the hybrid ACA algorithm \cite{Grasedyck2005HACA}. The hybrid ACA algorithm assumes prior knowledge about the skeleton rows and columns to leverage interpolation algorithms (e.g., ID and CUR) on a skeleton submatrix and use ACA to refine the skeletons. In contrast, BACA uses cross approximations with QRCP to select skeleton rows and columns and uses interpolation algorithms (LRID at line \ref{ID}) to form the low-rank update in each iteration. In other words, hybrid ACA can be treated as embedding ACA into interpolation algorithms while BACA can be thought of as embedding interpolation algorithms into ACA iterations. In addition, BACA is purely algebraic and requires no prior knowledge of the row/column skeletons or geometrical information about the rows/columns.}    
	
	It is worth mentioning that the choice of $d$ affects the trade-off between efficiency and robustness of the BACA algorithm. When $d<r$, the algorithm requires $O(nr^2)$ operations assuming convergence in $O(r/d)$ iterations as each iteration requires $O(nr_kd)$ operations. \ylrev{For example, BACA (Algorithm \ref{algo_baca}) precisely reduces to ACA (Algorithm \ref{algo_aca}) when $d=1$. In what follows we refer to the baseline ACA algorithm as BACA with $d=1$.} On the other hand, BACA converges in a constant number of iterations when $d\gg r$. \ylrev{In the extreme case, BACA reduces to QRCP-based ID when $d=\mathrm{min}\{m,n\}$ (note that the LRID algorithm at line \ref{ID} remains the only nontrivial operation)}. In this case the algorithm requires $O(n^2r)$ operations but enjoys the provable convergence of QRCP. \ylrev{Detailed complexity analysis of the BACA algorithm will be provided in Section \nameref{algo_ana}.} 
	
	The BACA algorithm oftentimes exhibits overestimated ranks compared to those revealed by truncated SVD. Therefore, an SVD re-compression step of $U$ and $V$ may be needed via first computing a QR of $U$ and $V$ as $[Q_{U},T_{U}]=\mathtt{QR}(U)$, $[Q_{V},T_{V}]=\mathtt{QR}(V^t)$, and then \ylrev{a truncated} SVD of $T_UT_V^t$ \cite{Heldring2015stochasticACA}. The result can be viewed as an approximate truncated SVD of $A$ and we assume this is the output of the BACA algorithm in the rest of this paper.    
	
	\IncMargin{1em}
	\begin{algorithm2e}
		\SetKwData{Left}{left}\SetKwData{This}{this}\SetKwData{Up}{up}\SetKwData{Conv}{conv}\SetKwData{Err}{er}
		\SetKwFunction{QR}{QR}\SetKwFunction{BACA}{BACA}\SetKwFunction{SVD}{SVD}\SetKwFunction{LRID}{LRID}\SetKwFunction{LRnorm}{LRnorm}\SetKwFunction{LRnormUp}{LRnormUp}
		\SetKwInOut{Input}{input}\SetKwInOut{Output}{output}
		\SetKwProg{Fn}{Function}{}{}
		\Input{Matrix $A\in\mathbb{R}^{m\times n}$, number of leaf-level subblocks $n_b$, block size $d$ of leaf-level BACA, relative tolerance $\epsilon$}
		\Output{\ylrev{Truncated SVD of $A\approx U\Sigma V$ with rank $r$}}
		Create $L$-level trees on index vectors $[1,m]$ and $[1,n]$ with index set $I_\tau$ and $J_\nu$ for nodes $\tau$ and $\nu$ at each level, $L=\mathrm{log}{\sqrt{n_b}}$, the leaf and root levels are denoted $0$ and $L$, respectively\;
		\For{$l= 0$ \KwTo $L$}{
			\ForEach{$A_{\tau\nu}=A(I_{\tau},J_{\nu}) \mathrm{~at~level~} l$}{
				\If{leaf-level}{
					$[U_{\tau\nu},\Sigma_{\tau\nu},V_{\tau\nu},r_{\tau\nu}]=\BACA(A_{\tau\nu},d,\epsilon)$\;
				}
				\Else{
					Let $\tau_1,\tau_2$ and $\nu_1,\nu_2$ denote children of $\tau$ and $\nu$\;
					\For{$i= 1$ \KwTo $2$}{
						$\bar{U}_{\tau_i\nu}$ = $[U_{\tau_i\nu_1}\Sigma_{\tau_i\nu_1},U_{\tau_i\nu_2}\Sigma_{\tau_i\nu_2}]$\;
						$\bar{V}_{\tau_i\nu}$ = $\diag(V_{\tau_i\nu_1},V_{\tau_i\nu_2})$\;
						$[U_{\tau_i\nu},\Sigma_{\tau_i\nu},{V}_{\tau_i\nu},r_{\tau_i\nu}]\leftarrow\SVD(\bar{U}_{\tau_i\nu},\epsilon)$\;
						$V_{\tau_i\nu}\leftarrow V_{\tau_i\nu}\bar{V}_{\tau_i\nu}$\;
					}
					$\bar{U}_{\tau\nu}$ = \ylrev{$\diag(U_{\tau_1\nu},U_{\tau_2\nu})$}\;
					$\bar{V}_{\tau\nu}$ = $[\Sigma_{\tau_1\nu}V_{\tau_1\nu};\Sigma_{\tau_2\nu}V_{\tau_2\nu}]$\;
					$[{U}_{\tau\nu},{\Sigma}_{\tau\nu},V_{\tau\nu},r_{\tau\nu}]\leftarrow\SVD(\bar{V}_{\tau\nu},\epsilon)$\;
					$U_{\tau\nu}\leftarrow \bar{U}_{\tau_\nu}U_{\tau\nu}$\;			
				}
			}
		}
		return $U={U}_{\tau\nu}$, $V={V}_{\tau\nu}$, $\Sigma={\Sigma}_{\tau\nu}$, $r=r_{\tau\nu}$\;
		\caption{Hierarchical low-rank merge algorithm with BACA (H-BACA)}\label{algo_merge}
	\end{algorithm2e}\DecMargin{1em}

	\subsection{Parallel Hierarchical Low-Rank Merge}
	The distributed-memory implementations of the proposed BACA algorithm and the baseline ACA algorithm can pose performance challenges as straightforward parallelization of all operations in Algorithm \ref{algo_baca} and \ref{algo_aca} involves many collective communications. To see this, assuming the $U$ and $V$ factors in Algorithm \ref{algo_aca} follow 1D block row and column data layouts, then every operation from line 3 to line 9 requires one or more collective communications. Instead, one can assign one process to perform BACA/ACA on submatrices without any communication and then leverage parallel low-rank arithmetic to merge the results into one single low-rank product. To elucidate the proposed algorithm, we first describe the hierarchical low-rank merge algorithm then outline its parallel implementation.     
	
\begin{figure*}[!tb]
		\centering	
		\begin{subfigure}{\subfigwidthww}
			\includegraphics[width=\subfigwidthww]{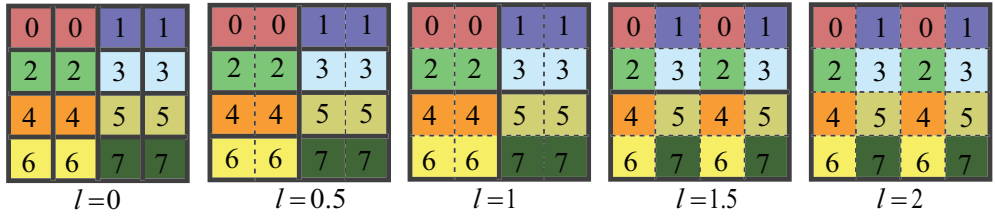}
			\vspace{-20pt}						
		\end{subfigure}	
		\caption{Parallel hierarchical merge with 8 processes. \ylrev{Blocks surrounded by solid lines represent $A_{\tau\nu}$ after compression at ech step $l$. Blocks surrounded by dashed lines represent ScaLAPACK blocks.}} \label{fig:parallel}	
	\end{figure*}

	Given a matrix $A\in\mathbb{R}^{m\times n}$ with $m\approx n$, the algorithm first creates $L$-level binary trees for index vectors $[1,m]$ and $[1,n]$ with index set $I_\tau$ and $J_\nu$ for nodes $\tau$ and $\nu$ at each level, \ylrev{upon recursively dividing each index set into $I_{\tau_i}$/$J_{\nu_j}$ of approximately equal sizes, $i=1,2$, $j=1,2$. Here, $\tau_i$ and $\nu_j$ are children of $\tau$ and $\nu$, respectively. The leaf and root levels are denoted $0$ and $L$, respectively.} This process generates $n_b$ leaf-level submatrices \ylrev{of similar sizes. For simplicity, it is assumed $n_b=4^{L}$.} We denote submatrices associated with $\tau,\nu$ as $A_{\tau\nu}=A(I_\tau,J_\nu)$ and their truncated SVD as $[U_{\tau\nu},\Sigma_{\tau\nu},V_{\tau\nu},r_{\tau\nu}]=\mathtt{SVD}(A_{\tau\nu},\epsilon)$. Here $r_{\tau\nu}$ is the $\epsilon$-rank of $A_{\tau\nu}$. \ylrev{As submatrices $A_{\tau\nu}$ have significantly smaller dimensions than $A$ (e.g., when $n_b=O(n^2)$ as an extreme case), both BACA and ACA algorithms become more robust to attain the truncated SVD. Following compression of $n_b$ submatrices $A_{\tau\nu}$ by BACA or ACA at step $l=0$, there are multiple approaches to combine them into one low-rank product including randomized algorithms via applying $A$ to random matrices, and deterministic algorithms via recursively pair-wise re-compressing the blocks using low-rank arithmetic. Here we choose the deterministic algorithm for simplicity of rank estimation and parallelization.} Here, we deploy truncated SVD as the re-compression tool but other tools such as ID, QR, UTV can also be applied. Fig. \ref{fig:merge} illustrates one re-compression operation for transforming SVDs of $A_{\tau_i\nu_j}, i=1,2, j=1,2$ into that of $A_{\tau\nu}$. The operation first horizontally compresses \ylrev{SVDs of} $A_{\tau_i\nu_j}, i=1,2, j=1,2$ \ylrev{at step $l-\frac{1}{2}$} and then vertically compresses the results, \ylrev{i.e., SVDs of} $A_{\tau_i\nu}, i=1,2$ \ylrev{at step $l$, $l=1,..,L$}. Specifically, the horizontal compression step is composed of one concatenation operation in (\ref{columncomb}) and one compression operation in (\ref{columnsvdmerge}):
	\begin{linenomath*} 
		\begin{align}
		\bar{U}_{\tau_i\nu} = [U_{\tau_i\nu_1}\Sigma_{\tau_i\nu_1},U_{\tau_i\nu_2}\Sigma_{\tau_i\nu_2}], ~\bar{V}_{\tau_i\nu} = \diag(V_{\tau_i\nu_1},V_{\tau_i\nu_2})\label{columncomb}\\
		[U_{\tau_i\nu},\Sigma_{\tau_i\nu},{V}_{\tau_i\nu},r_{\tau_i\nu}]\leftarrow\mathtt{SVD}(\bar{U}_{\tau_i\nu},\epsilon), ~V_{\tau_i\nu}\leftarrow V_{\tau_i\nu}\bar{V}_{\tau_i\nu}\label{columnsvdmerge}
		\end{align} 
	\end{linenomath*}

	with $i=1,2$. Let \ylrev{$\bar{U}_{\tau_i\nu}\bar{V}_{\tau_i\nu}$ and $U_{\tau_i\nu}\Sigma_{\tau_i\nu}V_{\tau_i\nu}$} denote the submatrix before \ylrev{and after} the SVD truncation, respectively. Similarly, the vertical compression step can be performed via horizontal merge of $A_{\tau_i\nu}^t, i=1,2$. \ylrev{Let $s_l$ represent the maximum rank $r_{\tau\nu}$ among all blocks at steps $l=0,1,...,L$.} Note that the algorithm returns an approximate truncated SVD after $L$ steps. \ylrev{As an example, the hierarchical merge algorithm with the level count of the hierarchical merge $L=2$ and $n_b=16$ is illustrated in Fig. \ref{fig:parallel}. At step $l=0$, the algorithm compresses all $n_b$ submatrices with BACA; at step $l=0.5,1.5$, the algorithm merges every horizontal pair of blocks; similarly at level $l=1,2$, the algorithm merges every vertical pair of blocks. Note that blocks surrounded by solid lines represent results after compression at each step $l$.} 

	The above-described hierarchical algorithm with BACA for leaf-level compressions, is dubbed H-BACA (Algorithm \ref{algo_merge}). In the following, a distributed-memory implementation of the H-BACA algorithm is described. Without loss of generality, it is assumed that $m=n=2^i$ \ylrev{and $p=2^j$}. The proposed parallel implementation first creates \ylrev{two $\lceil\mathrm{log}\sqrt{p}\rceil$}-level binary trees with $p$ denoting the total number of MPI processes. One process performs BACA compression of \ylrev{one or two} leaf-level submatrices and low-rank merge operations from the bottom up until it reaches a submatrix shared by more than one process. Then, all such blocks are handled by PBLAS and ScaLAPACK with BLACS process grids that aggregate those in corresponding submatrices. \ylrev{Consider the example in Fig. \ref{fig:parallel} with process count $p=8$.} The workload of each process is labeled with its process rank and highlighted with one color. The dashed lines represent the ScaLAPACK \ylrev{blocks}. First, \ylrev{BACA compressions and} merge operations at $l=0, 0.5$ are handled locally by one process without any communication. Next, merge operations at $l=1, 1.5, 2$ are handled by BLACS grids of $2\times1$, $2\times2$, and $4\times2$, respectively. For illustration purposes, we select the ScaLAPACK \ylrev{block} size in Fig. \ref{fig:parallel} as $n_0\times n_0$ where $n_0$ is the dimension of the finest-level submatrices in the hierarchical merge algorithm \ylrev{and $n=\sqrt{n_b}n_0$}. In this case, the only required data redistribution is from step $l=1$ to $l=1.5$. However, the ScaLAPACK \ylrev{block} size may be set to much smaller numbers in practice, requiring data redistribution at each row/column re-compression step. \ylrev{Similarly, the requirement of $m=n=2^i$ and $p=2^j$ is not needed in practice.}
	
	\section{Cost Analysis}\label{algo_ana}     
	\begin{table*}[!tp]
		\centering	
		\begin{tabular}{|c|c|c|}
			\hline
			& constant rank & increasing rank \\		
			& $s_l\ylrev{\approx} r$ & $s_l\ylrev{\approx}r/\sqrt{n_b}\times 2^l$\\
			\hline
			\ylrev{\rm{BACA} $d\leq s_0$}& $O(nr^2\sqrt{n_b})$ & \ylrev{$O(nr^2)/\sqrt{n_b}$}  
			\\
			\hline
			\ylrev{Merge compute}& \ylrev{$O(nr^2\sqrt{n_b})$} & \ylrev{$O(nr^2)$}  
			\\
			\hline
			\ylrev{Merge communicate}& $[O(r\mathrm{log}^2p),O(nr\mathrm{log}^2p/\sqrt{p})]$ & $[O(r\mathrm{log}p),O(nr\mathrm{log}p/\sqrt{p})]$  
			\\		
			\hline		
			\hline						
		\end{tabular}
		\caption{\ylrev{Flop counts and communication costs for the leaf-level compression and hierarchical merge operations in Algorithm \ref{algo_merge} for two classes of low-rank matrices}. $n$ and $r$ denote matrix dimension and rank. $d$ denotes the block size in BACA. $p$ and $n_b$ denote number of processes and leaf-level submatrices. $s_l$ denotes maximum ranks among all level-$l$ submatrices.}\label{tab:flop}
	\end{table*}

	\begin{table*}[!tp]
		\centering	
		\begin{tabular}{|c|c|c|c|c|}
			\hline
			Algorithm & ACA/$\rm{ACA}^{+}$ & Hyrbird-ACA & BACA& H-BACA\\		
			\hline
			Pivot count per iteration & 1 & 1 & $d$ & $n_bd$\\
			\hline
			Cost (constant rank) & $O(nr^2)$ & $O(nr^2)$ & $O(nr^2)$ & $O(nr^2\sqrt{n_b})$\\
			\hline
			Cost (increasing rank) & $O(nr^2)$ & $O(nr^2)$ & $O(nr^2)$ & $O(nr^2)$\\
			\hline
			Pre-selection of submatrices & no & yes & no & no\\
			\hline	
			\hline						
		\end{tabular}
		\caption{\ylrev{Comparisons between proposed BACA, H-BACA algorithms and existing ACA algorithms. Note that the algorithms show increasing robustness from left to right.}}\label{tab:acacompare}
	\end{table*}

	In this section, the costs for computation and communication of the proposed BACA and H-BACA algorithms are analyzed. 
	\subsection{Computational Cost}
	First, the costs for BACA can be summarized as follows. Assuming BACA converges in $O(\lceil {r}/{d} \rceil)$ iterations, each iteration performs entry evaluation from the residual matrices, QRCP for pivot selection, LRID for forming the LR product, and estimation of matrix norms. The entry evaluation computes $O(nd)$ entries each requiring $O(r_k)$ operations; QRCP on block rows requires $O(nd^2)$ operations; the LRID algorithm requires $O(ndd_k+d_kd^2)$ operations; norm estimation requires $O(nr_kd_k)$ operations. Summing up these costs, the overall cost for the BACA algorithm is 
	\begin{linenomath*} 
		\begin{align}
		c_{BACA} = \sum_{k=1}^{O(\lceil {r}/{d} \rceil)}(nd^2+nr_kd+d_kd^2) \nonumber\\
		\leq O(nd^2+rd^2+nrd)O(\lceil {r}/{d} \rceil)=O(nr^2)\label{costbaca}
		\end{align}
	\end{linenomath*} 	   
	Here we assume the block size $d\leq r$. Note that when $d\gg r$ (e.g., $d=O(n)$), it follows that the worst-case complexity is $c_{BACA} =O(n^2r)$ by bypassing the pivot selection step that causes the $nd^2$ term. In practice, one would always avoid the case of $d\gg r$.

	Next, the computational costs of the H-BACA algorithm are analyzed. The costs are analyzed for two cases of \ylrev{distributions of the maximum ranks $s_l$ at each level}, i.e., $s_l=r$ (ranks stay constant during the merge) and \ylrev{$s_l\approx 2^lr/\sqrt{n_b}=2^{l-L}r$} (rank increases by a factor of $2$ per level), \ylrev{$l=0,1,...,L$}. The constant-rank case is often valid for matrices with their numerical ranks independent of matrix dimensions \ylrev{(e.g., random low-rank matrices, matrices representing well-separated interactions from low-frequency and static wave equations and certain quantum chemistry matrices)}; the increasing-rank case holds true for matrices whose ranks depend polynomially (with order no bigger than 1) on the matrix dimensions (e.g., those arising from high-frequency wave equations, \ylrev{matrices representing near-field interactions from low-frequency and static wave equations, and certain classes of kernel methods on high dimensional data sets}). From the aforementioned
analysis of BACA, the computational costs for the leaf-level compression \ylrev{$c_b=c_{BACA}n_b$} are:
	\begin{linenomath*} 
		\begin{align}
		\ylrev{c_b} &= O\Big(\frac{n}{\sqrt{n_b}}s_0^2n_b\Big), ~~\mathrm{if}~\ylrev{d\leq s_0}\label{costaca}
		\end{align}  
	\end{linenomath*}  
	which represent the complexity with ACA when \ylrev{$n_b=1$}. 
	
	Let $n_l=2^ln/\sqrt{n_b}$ denote the size of submatrices $A_{\tau,\nu}$ at level $l$. The computational costs $c_m$ of hierarchical merge operations can be estimated as 
	\begin{linenomath*} 
		\begin{align}
		c_m = \sum_{l=1}^{L}O(4^{L-l}n_ls_l^2)\label{costmerge}
		\end{align}
	\end{linenomath*} 
	\ylrev{Accounting for the two cases of rank distributions, the computational costs for the leaf-level BACA and hierarchical merge operations of the H-BACA algorithm} are summarized in Table \ref{tab:flop}. \ylrev{Note that the costs of the BACA algorithm can also be extracted from Table \ref{tab:flop} upon setting $n_b=1$.} Not surprisingly, the hierarchical merge algorithm induces a computational overhead of \ylrev{at most} $\sqrt{n_b}$ when ranks stay constant; \ylrev{the leaf-level compression can have a $1/\sqrt{n_b}$ reduction factor for the increasing rank case and $\sqrt{n_b}$ overhead for the constant rank case.}
	
	\ylrev{For completeness, the comparison between the proposed BACA, H-BACA algorithms (assuming $d\leq r_0$) and existing ACA algorithms are given in Table \ref{tab:acacompare}. In contrast to existing ACA algorithms that select one pivot at a time, BACA and H-BACA select $d$ and $n_bd$ pivots simultaneously. As such, H-BACA is the most robust algorithm among all listed here. Not surprisingly, H-BACA can induce a computational overhead of $\sqrt{n_b}$.}       
	\subsection{Communication Cost}
	As the leaf-level BACA compression requires no communication, only the communication costs for the hierarchical merge operations are analyzed here. Since the merge operations may introduce an $O(\sqrt{n_b})$ computational overhead, one would only increase $n_b$ to create more parallelism, i.e., the process count $p\approx n_b$. Let $p_l=4^l$ denote the number of processes involved in one level $l$ merge operation, $l=1,...,L$. The operation requires redistribution between process grids of sizes $p_l$, $2p_l$ and $4p_l$ (see the example in Fig. \ref{fig:parallel}). Each process grid involves a PDGEMM function in PBLAS to combine the low-rank products and a PDGESVD function in ScaLAPACK to compute the new rank after the combination (see Fig. \ref{fig:merge}). Let the pair [\#messages, volume] denote the communication cost including the number of messages and the number of words transferred along the critical path. Then the communication costs for each (BLACS) grid redistribution, PDGEMM and PDGESVD during the hierarchical merge are $[O(1),O(n_ls_l/p_l)]$, $[O(s_l),O(n_ls_l/\sqrt{p_l})]$, and $[O(s_l\mathrm{log}p_l),O(n_ls_l\mathrm{log}p_l/\sqrt{p_l})]$, respectively. Recall that $n_l=2^ln/\sqrt{p}$ and $s_l$ denote the size and rank of submatrices at level $l$ and note that $n_l\gg s_l$. Therefore the communication cost $v_m$ of the hierarchical merge (and H-BACA) can be estimated as 
	\begin{linenomath*} 
		\begin{align}
		v_m = \sum_{l=1}^{L}\Big[O(s_l\mathrm{log}p_l),O\Big(\frac{n_ls_l\mathrm{log}p_l}{\sqrt{p_l}}\Big)\Big]\nonumber\\
		=\sum_{l=1}^{L}\Big[O(ls_l),O\Big(\frac{lns_l}{\sqrt{p}}\Big)\Big]\label{costvolume}
		\end{align}
	\end{linenomath*}      
	Consider the two cases of rank distributions, i.e., $s_l=r$ and $s_l\approx2^{l-L}r$, the overall communication costs of H-BACA are $v_m =[O(r\mathrm{log}^2p),O(nr\mathrm{log}^2p/\sqrt{p})]$ and $v_m =[O(r\mathrm{log}p),O(nr\mathrm{log}p/\sqrt{p})]$, respectively (see Table \ref{tab:flop}). 
	
	\section{Numerical Results}\label{numerical}
	This section presents several numerical results to demonstrate the accuracy and efficiency of the proposed H-BACA algorithm. The matrices in all numerical examples are generated from the following kernels: 1. Gaussian kernel: $A_{i,j}=\exp(\frac{-\left\lVert x_i-x_j\right\rVert^2}{2h^2})$, $i,j=1,...,2n$. Here $h$ is the Gaussian width, and $x_i\in\mathbb{R}^{8\times1}$ and $\mathbb{R}^{784\times1}$ are feature vectors in one subset of the SUSY and MNIST Data Sets from the UCI Machine Learning Repository \cite{Dua2017UCI}, respectively. Note that the Gaussian kernel permits low-rank compression as shown in \cite{RuoxiWang2017RBFrank,Bach2013Analysis,Musco2017Nystrom} 2. EFIE2D kernel: $A_{i,j}=H_0^{(2)}(k\left\lVert x_i-x_j\right\rVert)$ resulting from the Nystr\"{o}m discretization of the electric field integral equation (EFIE) for electromagnetic scattering from 2-D curves. Here $H_0^{(2)}$ is the second kind Hankel function of order 0, $k$ is the free-space wavenumber, $x_i, x_j\in\mathbb{R}^{2\times1}$ are discretization points (15 points per wavelength) of two 2-D parallel strips of length $1$ and distance $1$. 3. EFIE3D kernel: $A$ is obtained by the Galerkin method for EFIE to analyze electromagnetic scattering from 3-D surfaces. 4. Frontal3D kernel: $A$ is a dense frontal matrix that arises from the multifrontal sparse elimination for the finite-difference frequency-domain solution of the homogeneous-coefficient Helmholtz equation inside a unit cube. 5. Polynomial kernel: $A_{i,j}=(x_i^tx_j+h)^2$. Here $x_i, x_j\in\mathbb{R}^{50\times1}$ are points from a randomly generated dataset, and $h$ is a regularization parameter. \ylrev{6. Product-of-random kernel: $A=UV$ with $U\in\mathbb{R}^{n\times r}$ and $V\in\mathbb{R}^{r\times n}$ being random matrices with i.i.d. entries.} Note that the EFIE2D, EFIE3D and Frontal3D kernels result in complex-valued matrices. Throughout this section, we refer to ACA as a special case of BACA when $d=1$. In all examples \ylrev{except for the Product-of-random kernel,} the algorithm is applied to the offdiagonal submatrix $\ylrev{A_{12}}=A(1:n,1+n:2n)$ assuming rows/columns of $A$ have been properly permuted (e.g., by a KD-tree partitioning scheme). \ylrev{Note that the permutation may yield a hierarchical matrix representation of $A$, but in this paper we only focus on compression of one off-diagonal subblock of $A$ with H-BACA.} All experiments are performed on the Cori Haswell machine at NERSC, which is a Cray XC40 system and consists of 2388 dual-socket nodes with Intel Xeon E5-2698v3 processors running 16 cores per socket. The nodes are configured with 128\GB of DDR4 memory at 2133\MHz.

	\begin{figure*}[!tb]
		\centering	
		\begin{subfigure}{\subfigwidthw}
			\includegraphics[width=\subfigwidthw]{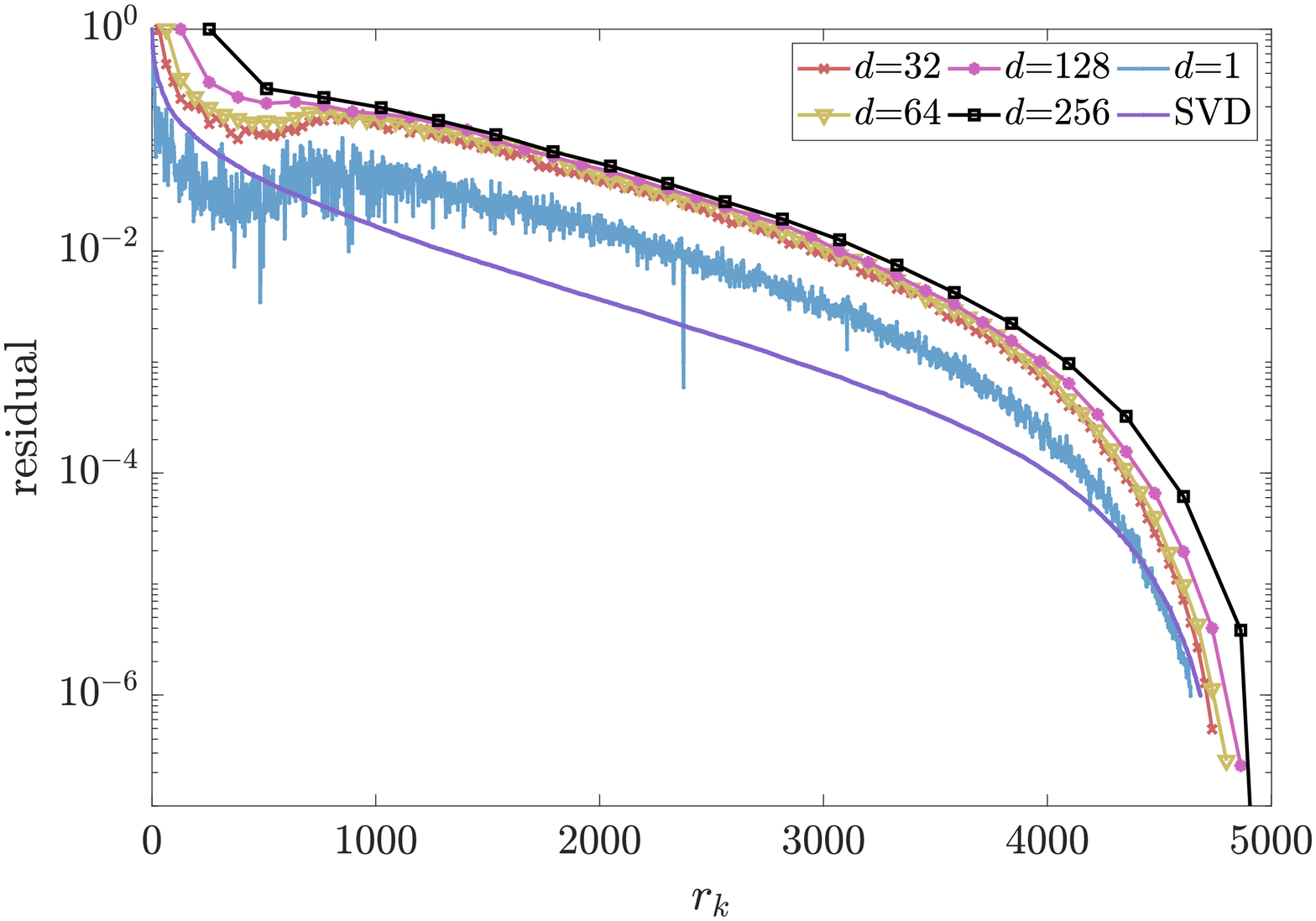}
			\vspace{-15pt}			
			\caption{Gaussian-SUSY ($h=1.0$)}\label{fig:history:gauss1}				
		\end{subfigure}
		\begin{subfigure}{\subfigwidthw}
			\includegraphics[width=\subfigwidthw]{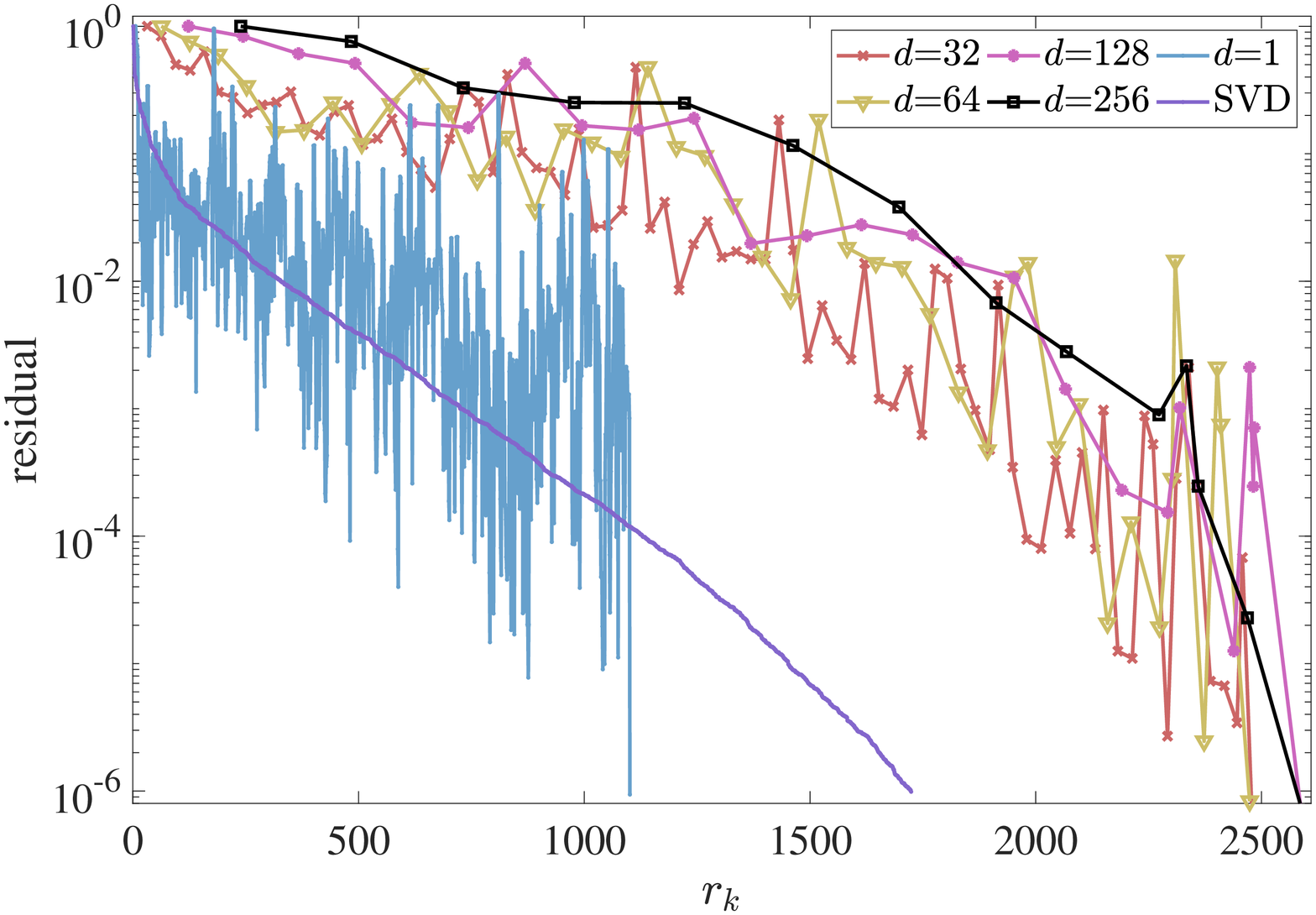}
			\vspace{-15pt}			
			\caption{Gaussian-SUSY( ($h=0.2$)}\label{fig:history:gauss02}		
		\end{subfigure}	
		\begin{subfigure}{\subfigwidthw}
			\includegraphics[width=\subfigwidthw]{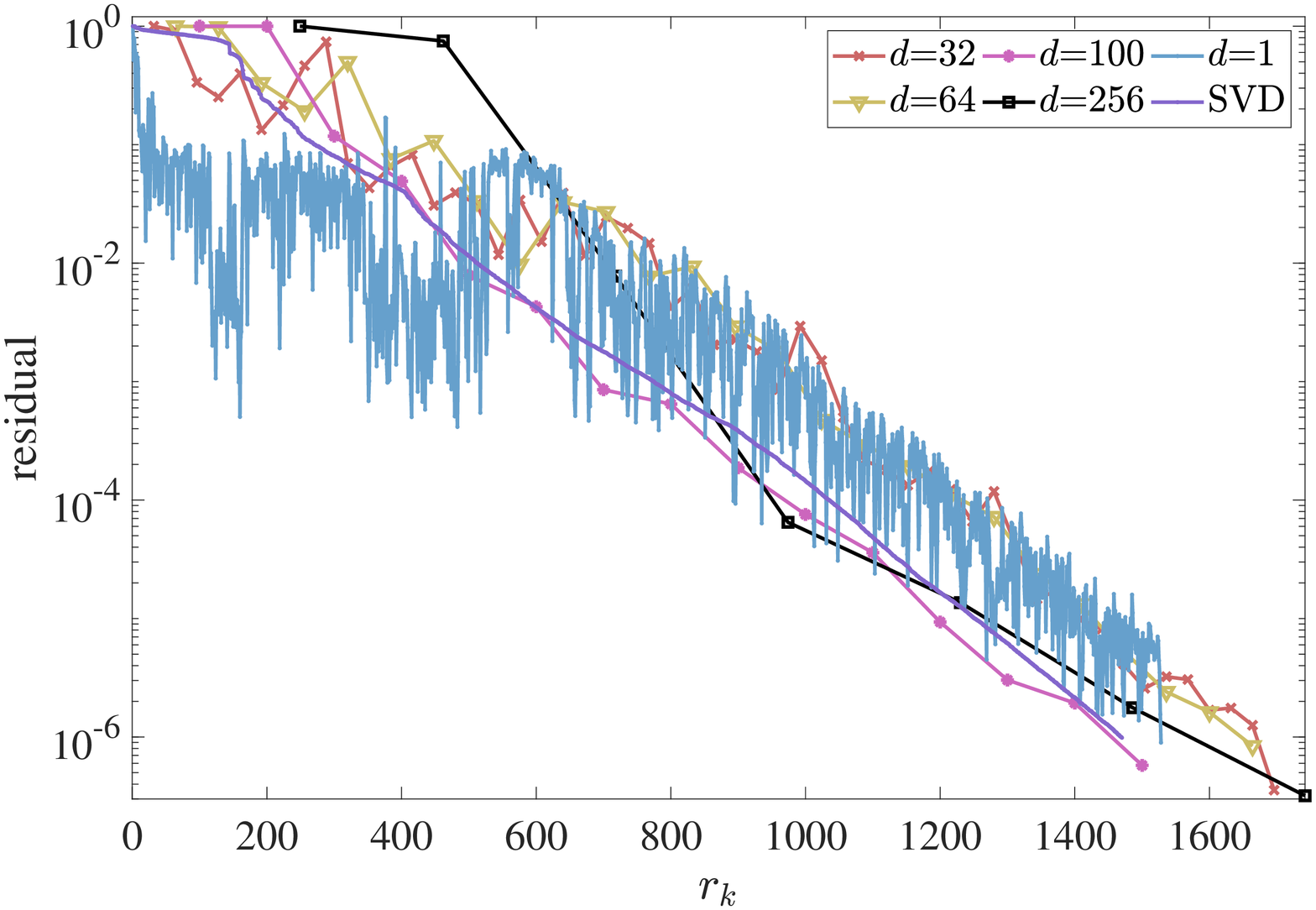}
			\vspace{-15pt}			
			\caption{EFIE3D}\label{fig:history:em3d}				
		\end{subfigure}
		\begin{subfigure}{\subfigwidthw}
			\includegraphics[width=\subfigwidthw]{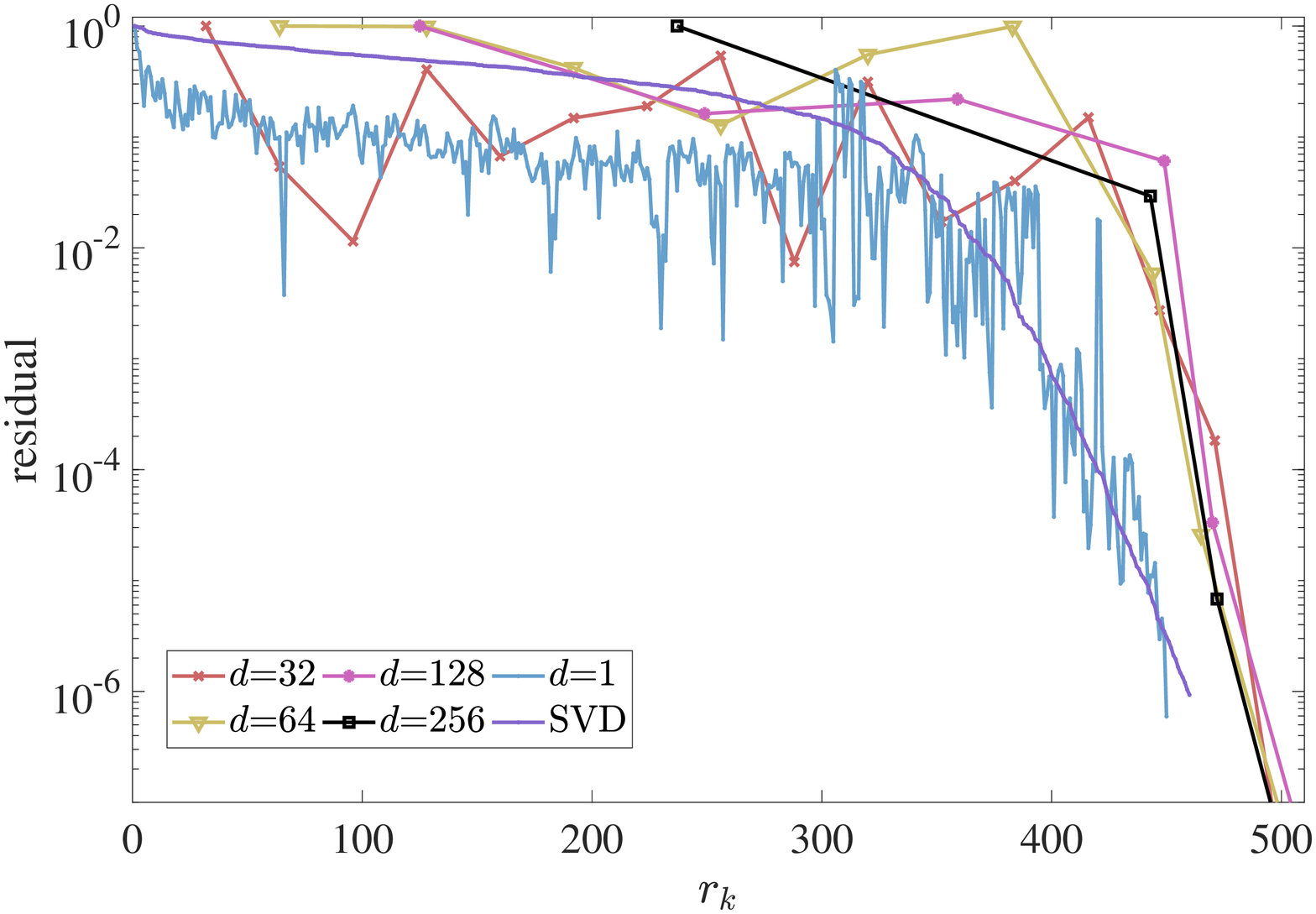}
			\vspace{-15pt}			
			\caption{Frontal3D}\label{fig:history:frontal}				
		\end{subfigure}
		\vspace{-5pt}	
		\caption{Convergence history of BACA for the (a) Gaussian-SUSY kernel with $h=1.0$, $n=5000$, $\epsilon=10^{-6}$, $r=4683$, (b) Gaussian-SUSY kernel with $h=0.2$, $n=5000$, $\epsilon=10^{-6}$, $r=1723$, (c) EFIE3D kernel for a unit sphere with $n=21788$, $\epsilon=10^{-6}$, $r=1488$ and (d) Frontal3D kernel with $n=1250$, $\epsilon=10^{-6}$, $r=718$} \label{fig:history}
	\end{figure*}

	\begin{figure*}[!tb]
		\centering	
		\begin{subfigure}{\subfigwidthw}
			\includegraphics[width=\subfigwidthw]{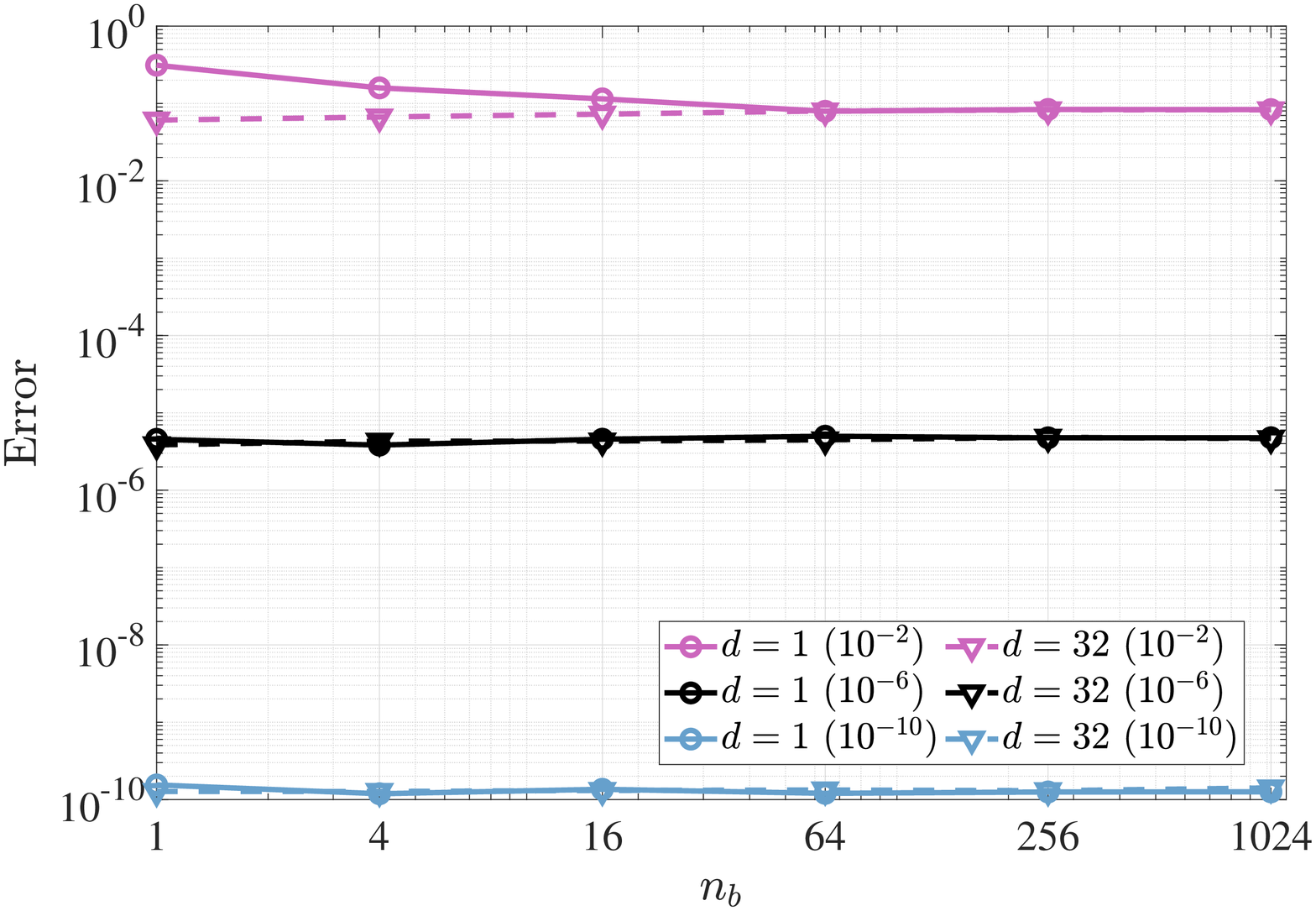}
			\vspace{-15pt}			
			\caption{Gaussian-SUSY ($h=1.0$)}\label{fig:accmerge:a}				
		\end{subfigure}
		\begin{subfigure}{\subfigwidthw}
			\includegraphics[width=\subfigwidthw]{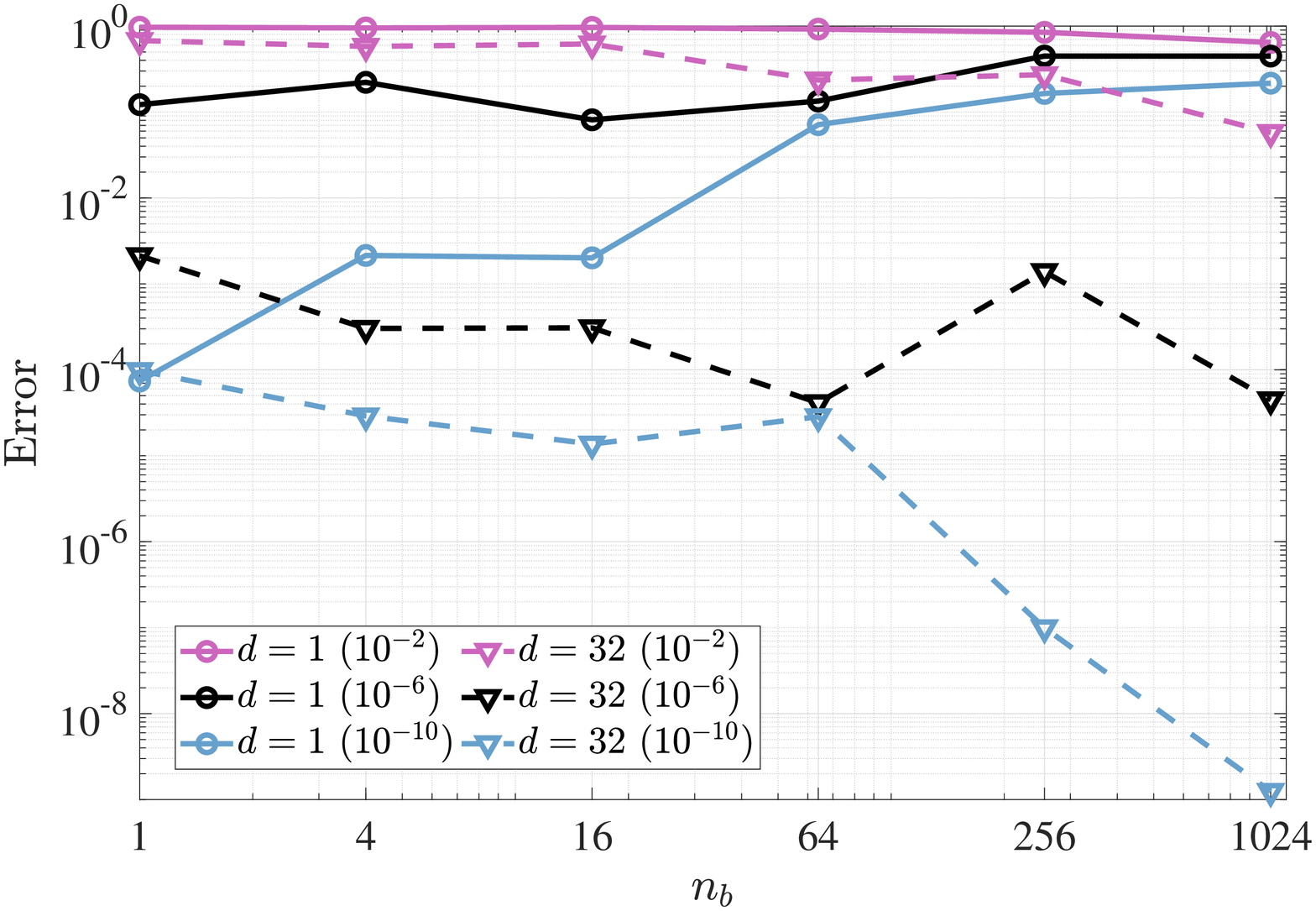}
			\vspace{-15pt}			
			\caption{Gaussian-SUSY ($h=0.2$)}\label{fig:accmerge:b}		
		\end{subfigure}	
		\begin{subfigure}{\subfigwidthw}
			\includegraphics[width=\subfigwidthw]{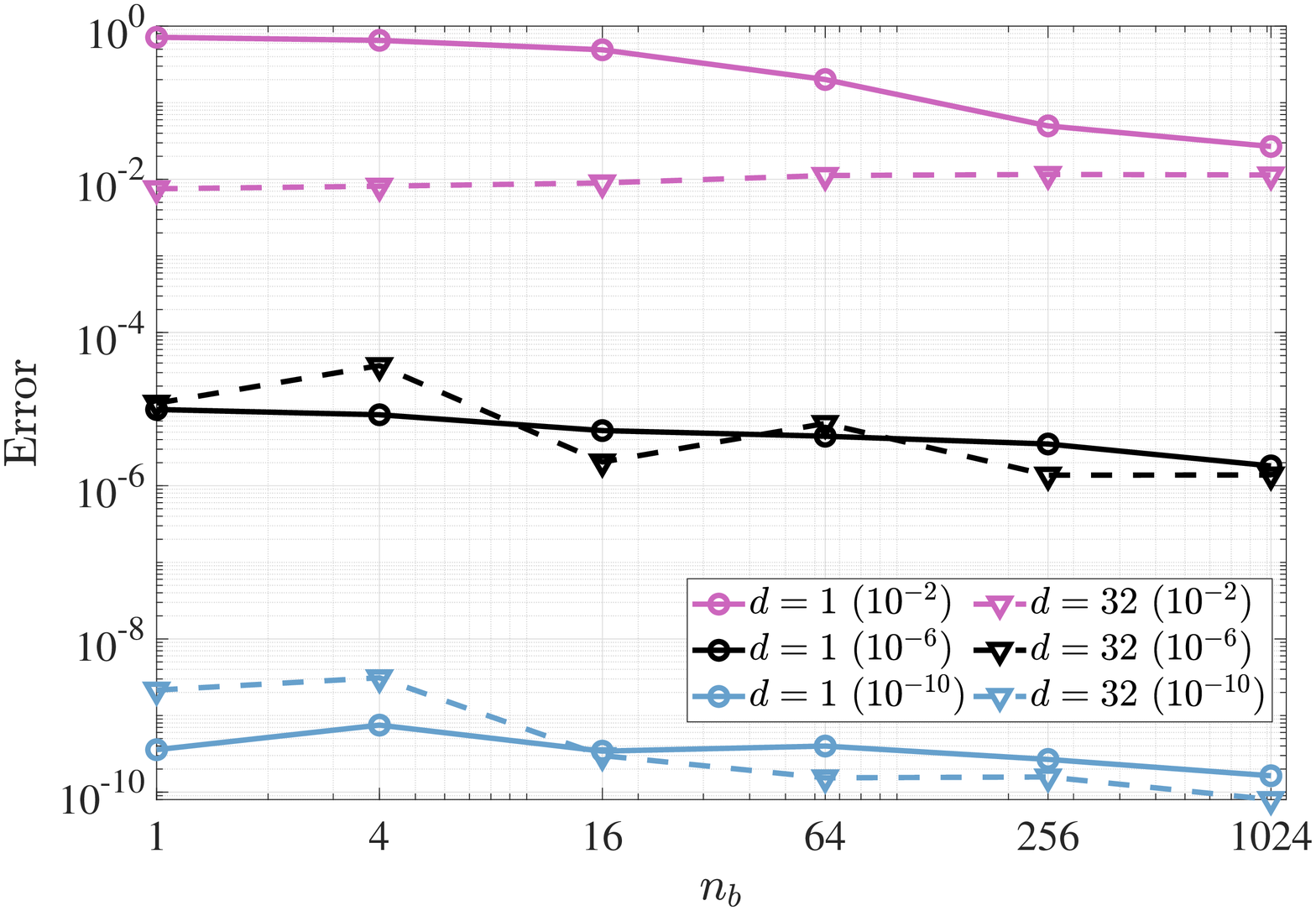}
			\vspace{-15pt}			
			\caption{EFIE3D}\label{fig:accmerge:c}		
		\end{subfigure}	
		\begin{subfigure}{\subfigwidthw}
			\includegraphics[width=\subfigwidthw]{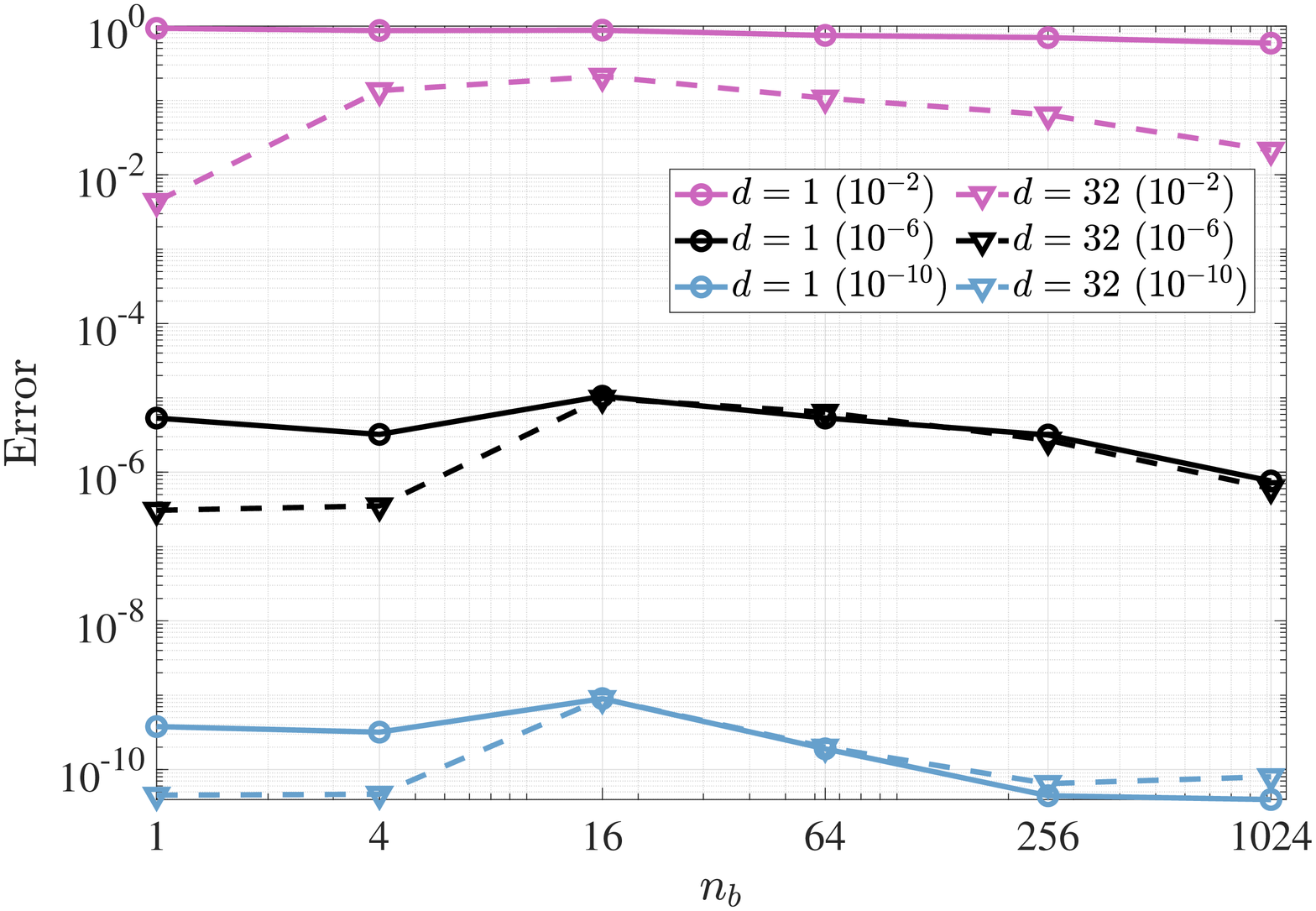}
			\vspace{-15pt}			
			\caption{Frontal3D}\label{fig:accmerge:d}		
		\end{subfigure}	
		\vspace{-5pt}	
		\caption{Measured error of H-BACA with $\epsilon=10^{-2},10^{-6},10^{-10}$ for the (a) Gaussian-SUSY kernel with $h=1.0$, $n=5000$, (b) Gaussian-SUSY kernel with $h=0.2$, $n=5000$ (c) EFIE3D kernel for a unit sphere with $n=1707$ and (d) Frontal3D kernel with $n=1250$.} \label{fig:accmerge}
	\end{figure*}

	\begin{figure*}[!tb]
		\centering	
		\begin{subfigure}{\subfigwidthw}
			\includegraphics[width=\subfigwidthw]{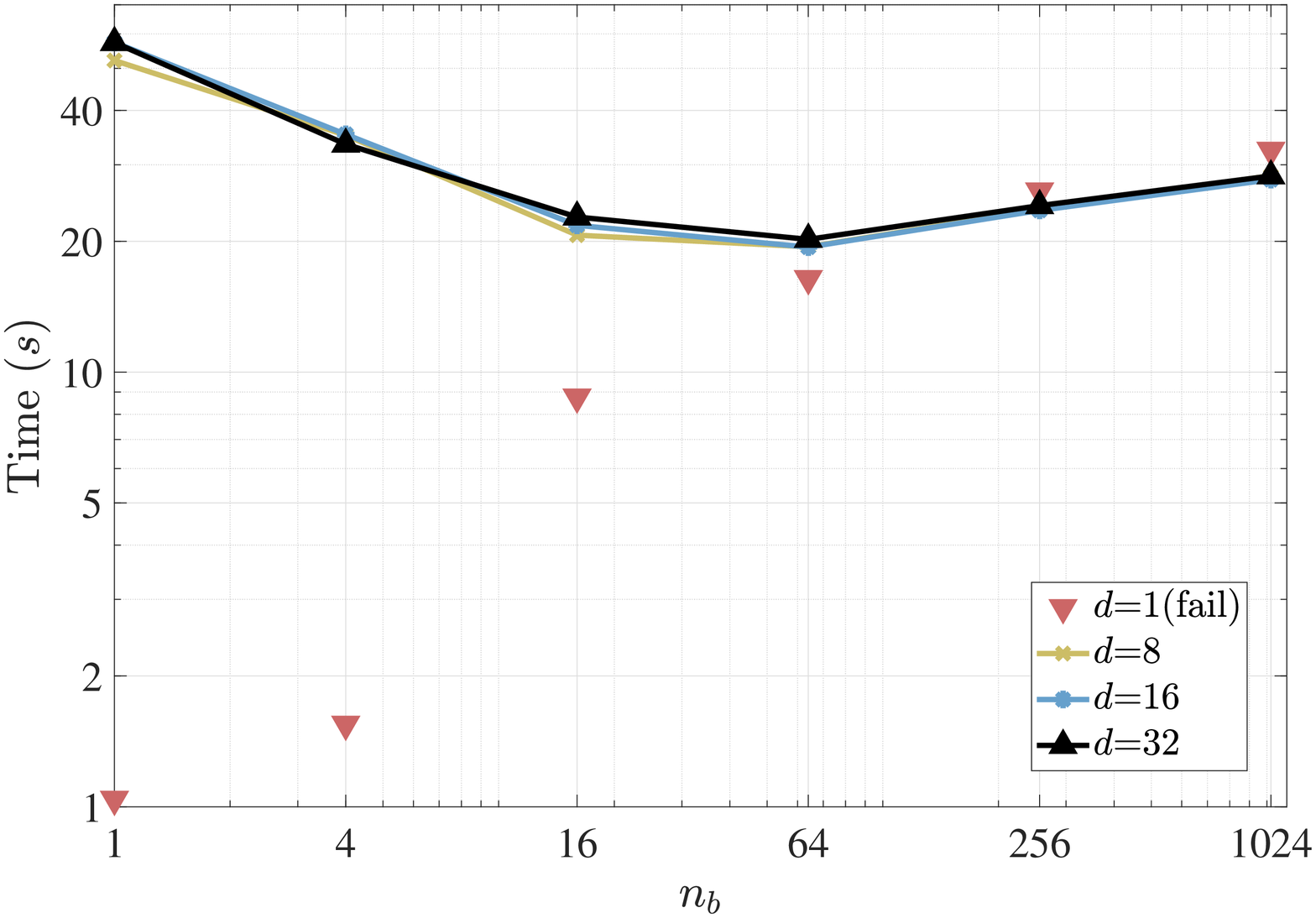}
			\vspace{-15pt}			
			\caption{\ylrev{Gaussian-SUSY}}\label{fig:many:a}		
		\end{subfigure}	
		\begin{subfigure}{\subfigwidthw}
			\includegraphics[width=\subfigwidthw]{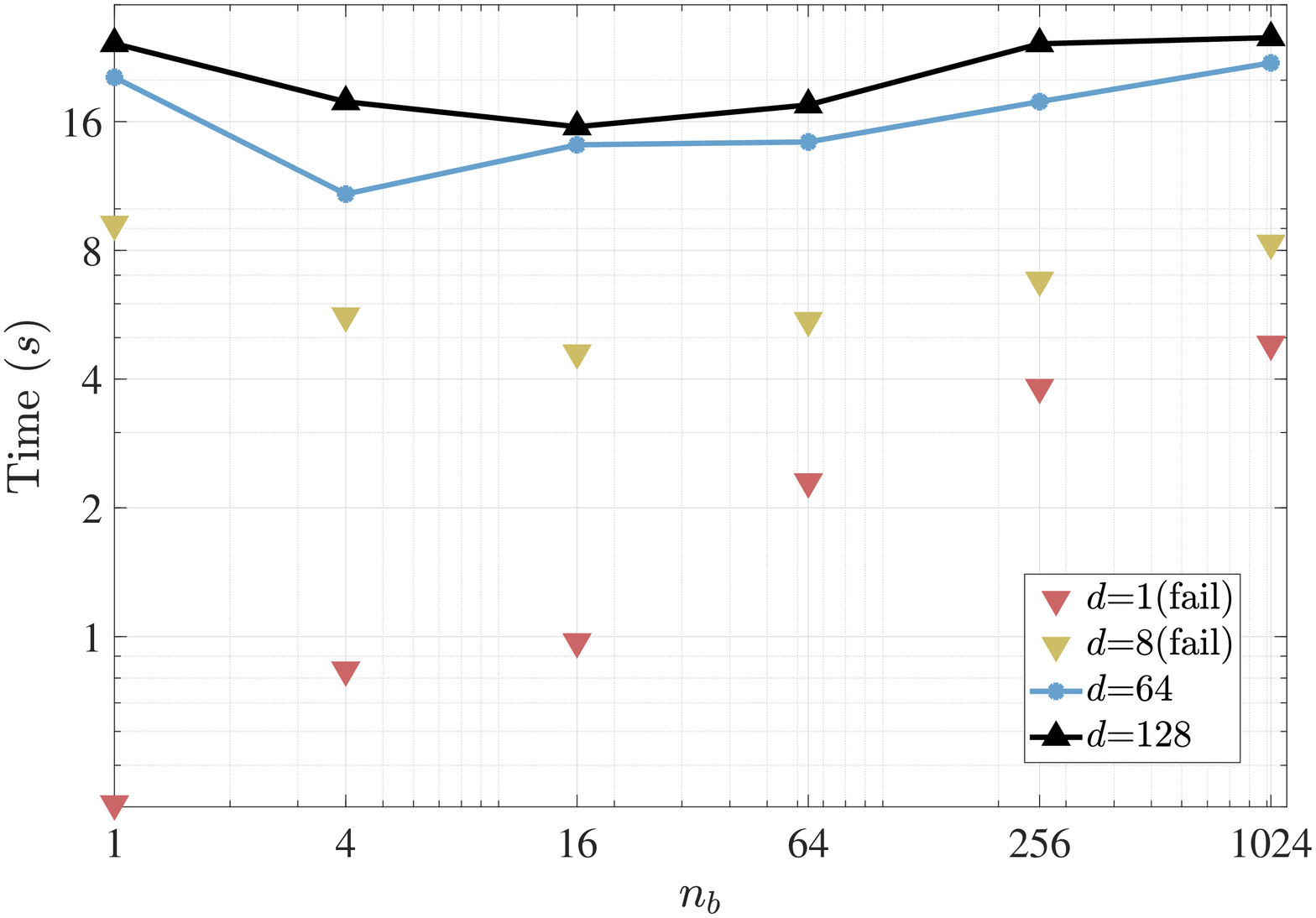}
			\vspace{-15pt}			
			\caption{\ylrev{Gaussian-MNIST}}\label{fig:many:b}		
		\end{subfigure}
		\begin{subfigure}{\subfigwidthw}
			\includegraphics[width=\subfigwidthw]{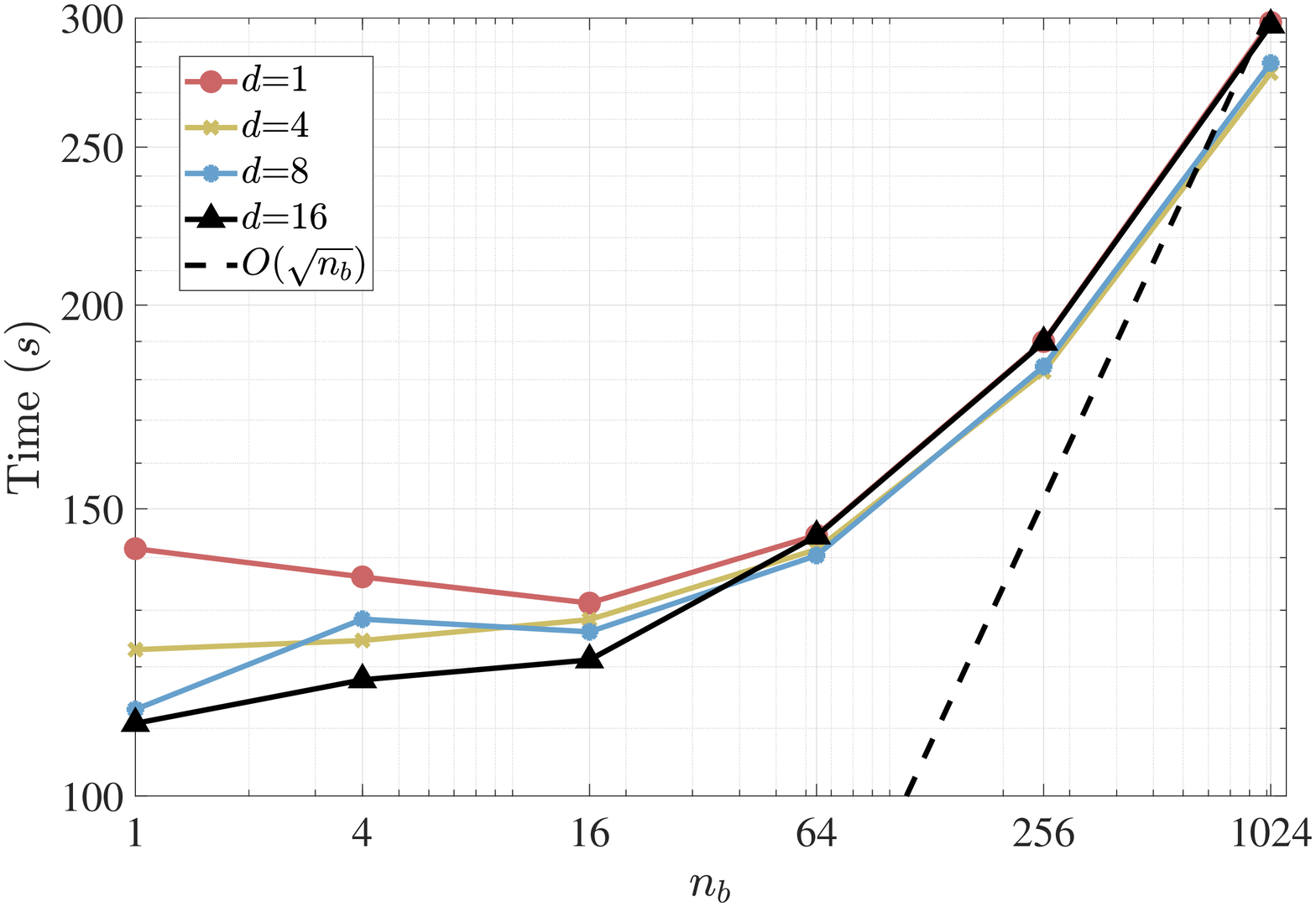}
			\vspace{-15pt}			
			\caption{\ylrev{EFIE3D}}\label{fig:many:c}		
		\end{subfigure}
		\begin{subfigure}{\subfigwidthw}
			\includegraphics[width=\subfigwidthw]{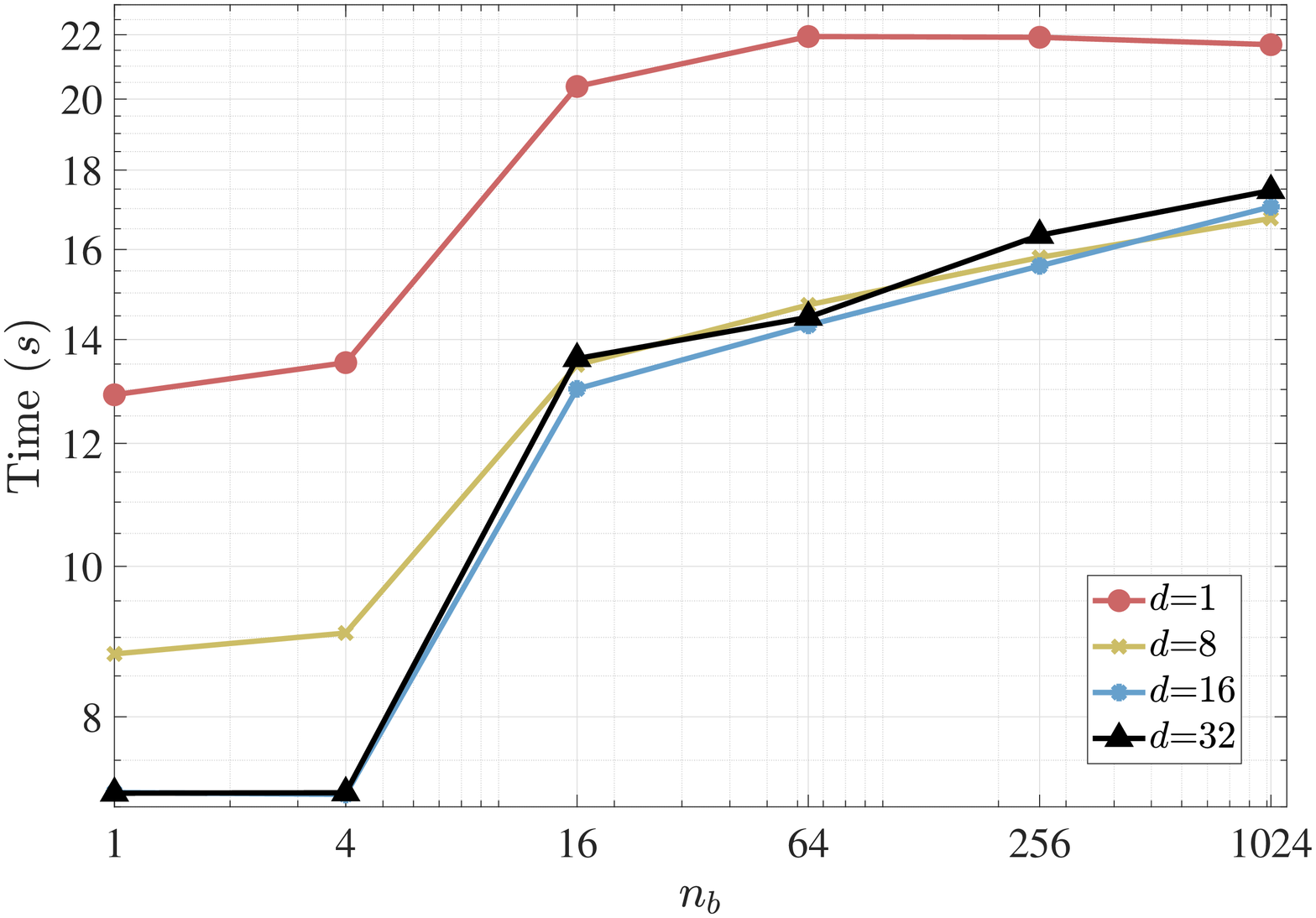}
			\vspace{-15pt}			
			\caption{\ylrev{Frontal3D}}\label{fig:many:d}		
		\end{subfigure}
		\begin{subfigure}{\subfigwidthw}
			\includegraphics[width=\subfigwidthw]{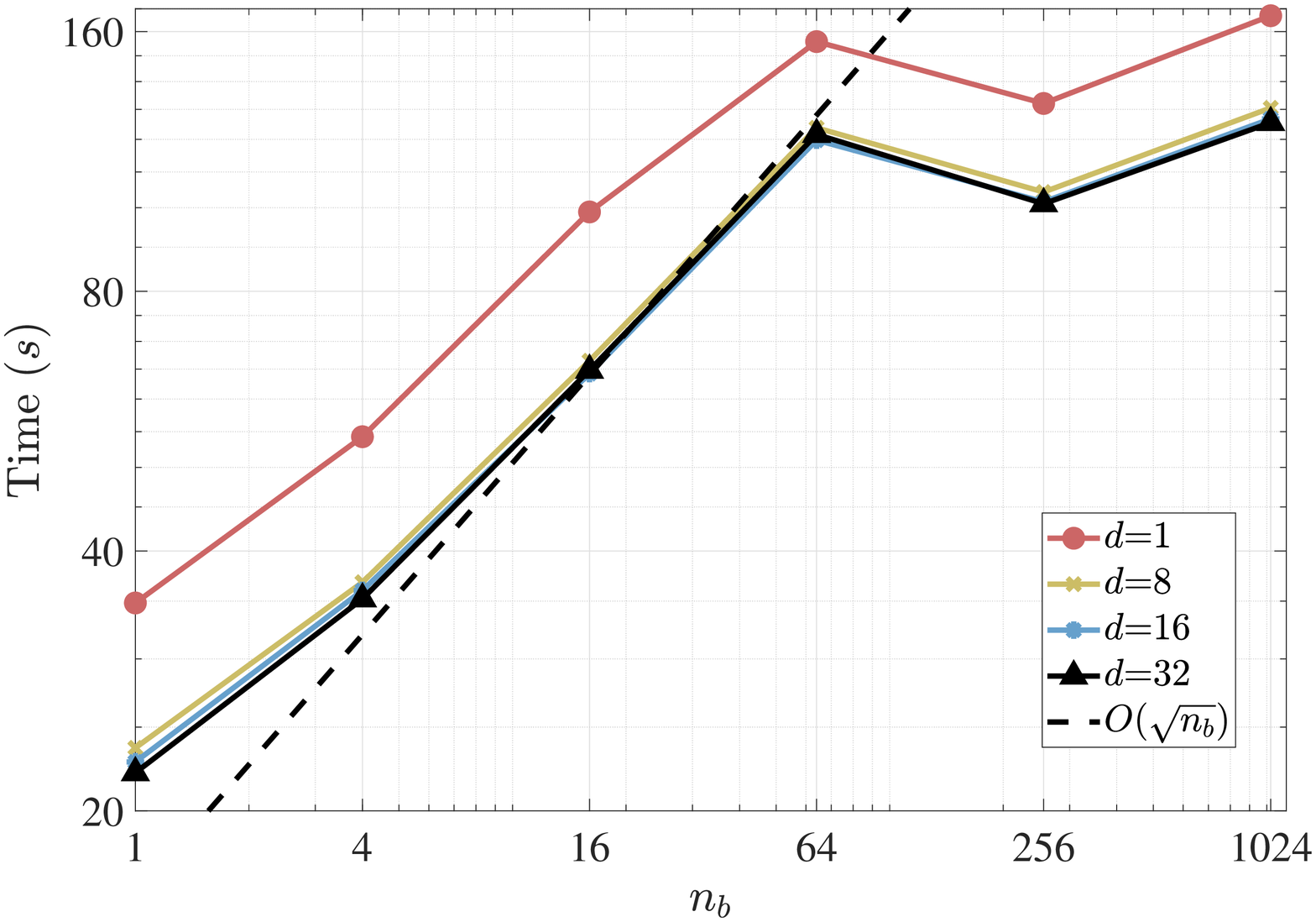}
			\vspace{-15pt}			
			\caption{\ylrev{Polynomial}}\label{fig:many:e}		
		\end{subfigure}
		\begin{subfigure}{\subfigwidthw}
			\includegraphics[width=\subfigwidthw]{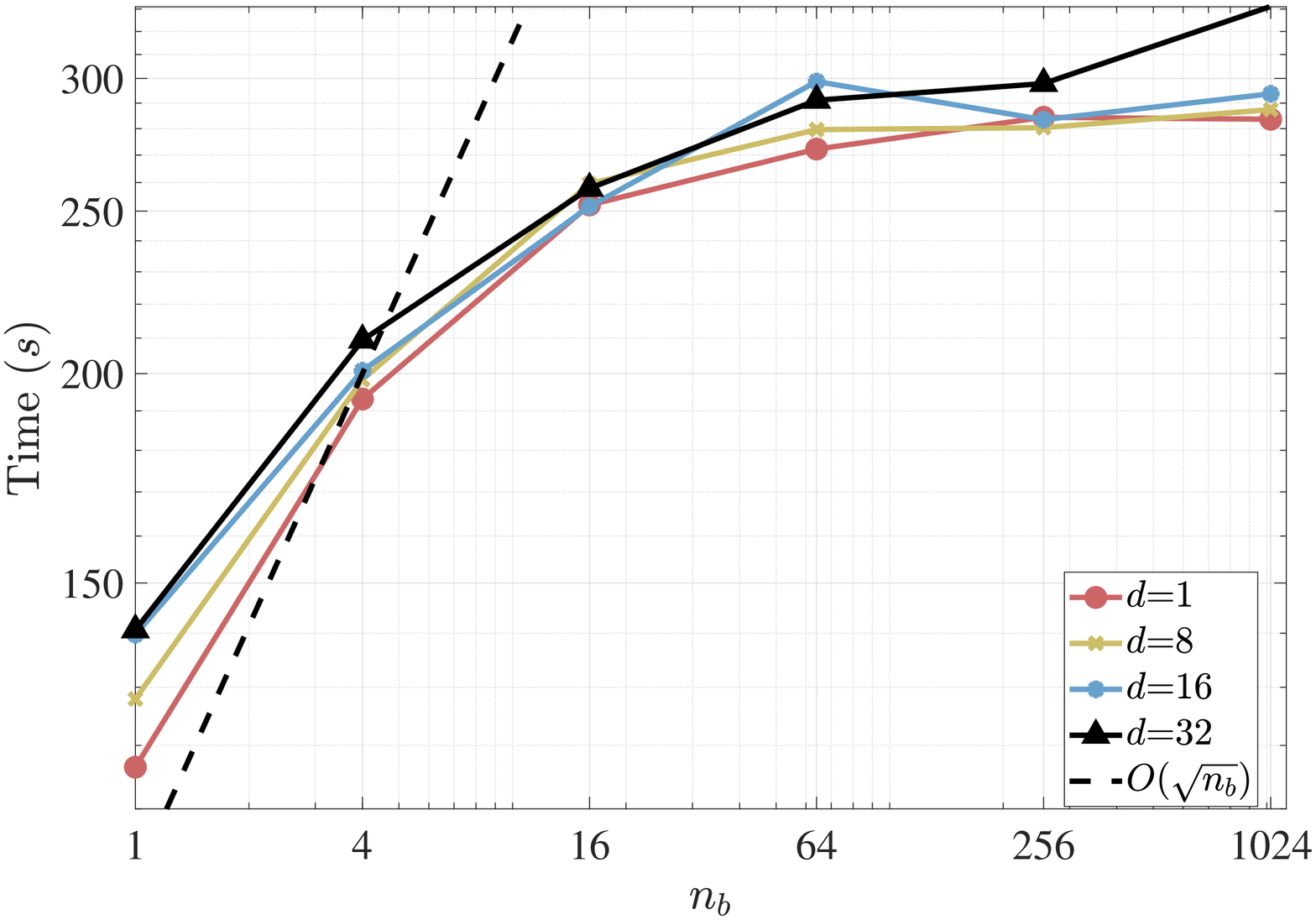}
			\vspace{-15pt}			
			\caption{\ylrev{Product-of-random}}\label{fig:many:f}		
		\end{subfigure}
		
		\vspace{-5pt}	
		\caption{Computation time of H-BACA with varying $n_b$ and $d$ for the (a) Gaussian-SUSY kernel with $h=1.0$, $n=50000$, $\epsilon=10^{-2}$, $r=298$, (b) Gaussian-MNIST kernel with $h=3.0$, $n=5000$, $\epsilon=10^{-2}$, $r=137$, (c) EFIE3D kernel for a unit sphere with $n=26268$, $\epsilon=10^{-6}$, $r=1488$, (d) Frontal3D kernel with $n=1250$, $\epsilon=10^{-6}$, $r=788$, (e) Polynomial kernel with $h=0.2$, $n=10000$, $\epsilon=10^{-4}$, $r=450$, and (f) Product-of-random kernel with $n=2500$, $r=1000$. Note that the data points where the algorithm fails \ylrev{are shown as triangular markers without lines}.}\label{fig:many}
	\end{figure*}

	\begin{figure*}[!tb]
		\centering	
		\begin{subfigure}{\subfigwidthw}
			\includegraphics[width=\subfigwidthw]{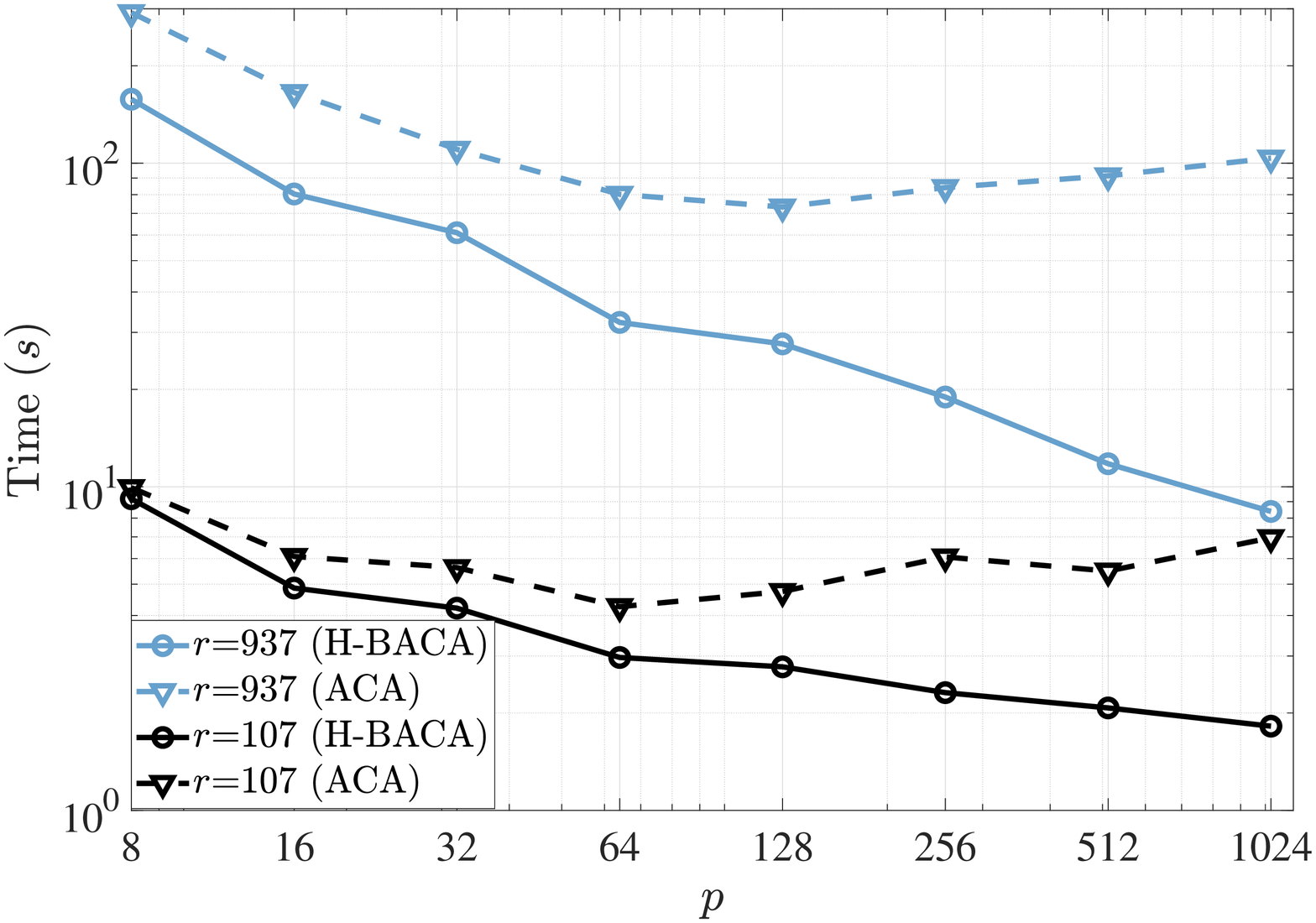}
			\vspace{-15pt}			
			\caption{EFIE2D}\label{fig:ss:a}				
		\end{subfigure}
		\begin{subfigure}{\subfigwidthw}
			\includegraphics[width=\subfigwidthw]{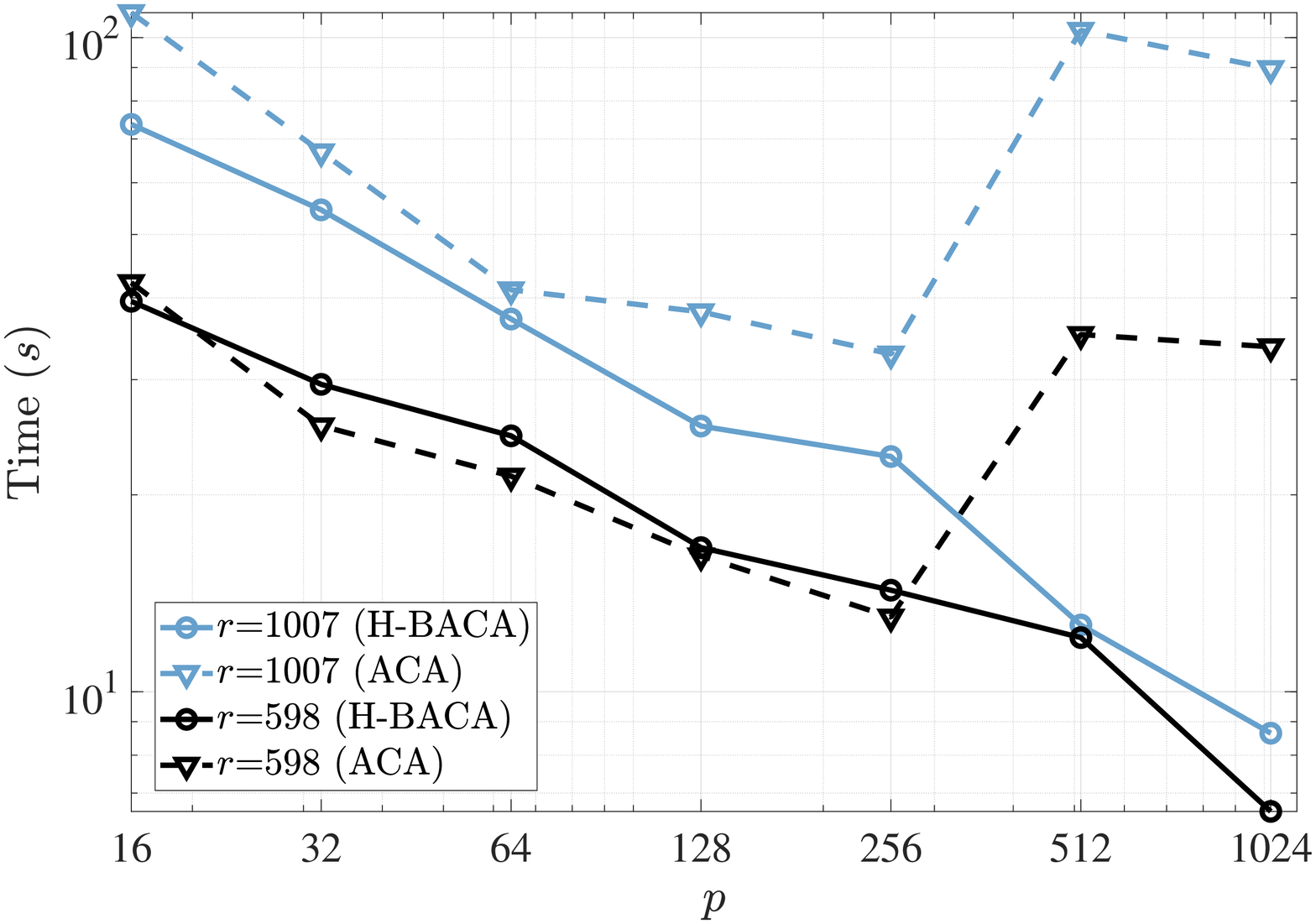}
			\vspace{-15pt}			
			\caption{\ylrev{EFIE3D}}\label{fig:ss:b}		
		\end{subfigure}
		\begin{subfigure}{\subfigwidthw}
			\includegraphics[width=\subfigwidthw]{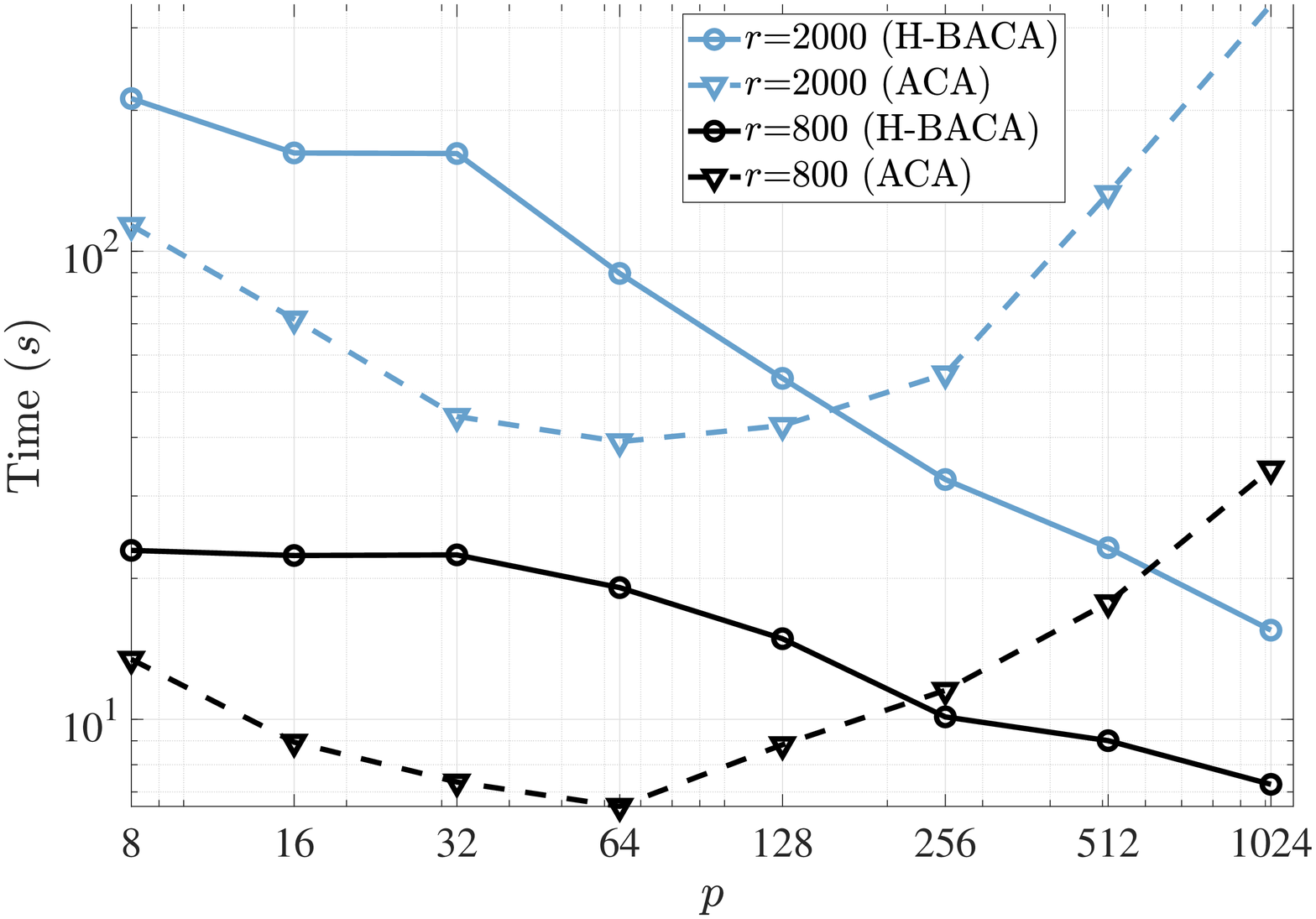}
			\vspace{-15pt}			
			\caption{\ylrev{Product-of-random}}\label{fig:ss:c}		
		\end{subfigure}
		\begin{subfigure}{\subfigwidthw}
			\includegraphics[width=\subfigwidthw]{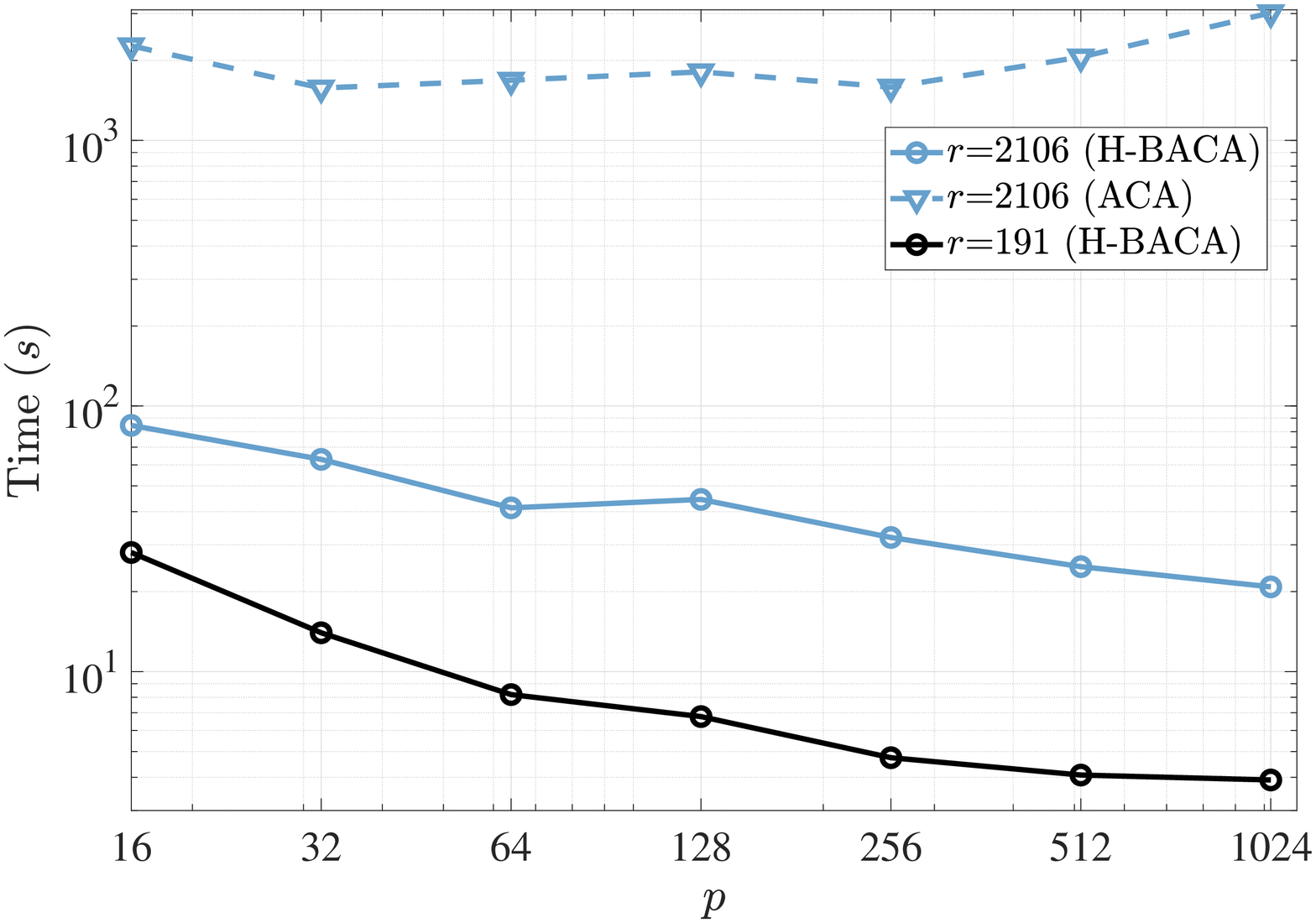}
			\vspace{-15pt}			
			\caption{\ylrev{Gaussian}}\label{fig:ss:d}		
		\end{subfigure}			
		\vspace{-5pt}	
		\caption{\ylrev{Computation time of H-BACA with varying process counts for the (a) EFIE2D kernel with $n=160000$, $\epsilon=10^{-4}$, $r=107, 937$ (b) EFIE3D kernel for a unit square with $n=21788$, $\epsilon=10^{-6}$, $r=598, 1007$, (c) Product-of-random kernel with $n=10000$, $r=800,2000$, and (d) Gaussian kernel for a randomly generated dataset with $h=1.0, 1.6$, $\epsilon=10^{-3}$, $r=2106, 191$. Note that for the Guassian matrix with $r=191$, ACA fails to provide accurate results and is not plotted. }}\label{fig:ss}
	\end{figure*}

	\subsection{Convergence}
	%\xslnote{Is this Section really about BACA algorithm, not H-BACA?}\ylnote{Yes. All the experiments in Fig 3 uses $n_b=1$. We can potentially also show the iteration history of BACA of all leaf-level submatrices in H-BACA. But I'm not sure how to make fair comparisons between different $n_b$}
	First, the convergence of the proposed BACA algorithm is investigated using several matrices: Gaussian-SUSY matrices with $n=5000$, $h=1.0, 0.2$, an EFIE3D matrix for a unit sphere with $n=21788$ and approximately 20 points per wavelength, and a Frontal3D matrix with $n=1250$ and 10 points per wavelength. The corresponding $\epsilon$-ranks are $r=4683, 1723, 1488, 718$ for $\epsilon=10^{-6}$. The residual histories versus revealed ranks $r_k$, at each iteration $k$ of BACA with $1\leq d\leq256$ are plotted in Fig. \ref{fig:history}. The residual error is defined as $\left\lVert U_kV_k\right\rVert_F/\left\lVert UV\right\rVert_F$ from (\ref{numu}). As a reference, the singular value spectra $\Sigma(k,k)/\Sigma(1,1)$ computed from $[U,\Sigma,V,r]=\mathtt{SVD}(A,\epsilon)$ are also plotted.

	For the Gaussian-SUSY matrices, the baseline ACA algorithm ($d=1$) behaves poorly
	with smaller $h$ due to the exponential decay of the Gaussian kernel. \ylrev{As a result, the matrix becomes increasingly sparse and coherent for small $h$ particularly for high dimensional data sets.} In fact, ACA constantly selects smaller pivots and the residual exhibits wild oscillations particularly for smaller $h$ (e.g., when $h=0.2$ in Fig. \ref{fig:history:gauss02}). Similarly, the analytical and numerical Green's functions respectively for the EFIE3D (Fig. \ref{fig:history:em3d}) and Frontal3D (Fig. \ref{fig:history:frontal}) matrices are not asymptotically smooth for ACA to converge rapidly. For all examples in Fig. \ref{fig:history}, significant portions of the residual curves lie below the singular value spectra which causes premature iteration termination for certain given residual errors. In stark contrast, the proposed BACA algorithm ($d=32,64,100,128,256$) shows increasingly smooth residual histories residing above the singular value spectra as the block size $d$ increases. Although BACA may overestimate the matrix ranks particularly for larger $d$, the SVD re-compression step mentioned in Section \nameref{baca_alg} can effectively reduce the ranks. 
	
	\subsection{Accuracy}
	Next, the accuracy of the H-BACA algorithm is demonstrated using the following matrices: two Gaussian-SUSY matrices with $n=5000$, $h=1.0,0.2$, one EFIE3D matrix for a unit sphere with $n=1707$ and approximately 20 points per wavelength, and a Frontal3D matrix with $n=1250$ and 10 points per wavelength. The relative Frobenious-norm error $\left\lVert A-UV\right\rVert_F/\left\lVert A\right\rVert_F$ is computed for changing number of leaf-level submatrices $n_b$ and block size $d$. When $h=1.0$ for the Gaussian-SUSY matrix (Fig. \ref{fig:accmerge}a), the H-BACA algorithms achieve desired accuracies ($\epsilon=10^{-2},10^{-6},10^{-10}$) using
	the baseline ACA ($d=1$), and BACA ($d=32$) when $n_b=1$
	%\xslnote{This seems to say that at the leaf level, H-BACA can use 3 different 
	%algorithms: ACA, BACA, QRCP?  But Algorithm 3 shows only BACA at the leaf level.}\ylnote{I added a sentence in first paragraph of section 5 mentioning ACA and QRCP are special cases of BACA. This was also discussed in the paragraph below equation 12}
	and the hierarchical merge operation only causes slight error increases as $n_b$ increases. However when $h=0.2$ for the Gaussian-SUSY matrix (Fig. \ref{fig:accmerge}b), all data points for H-BACA with $d=1$ fail due to the wildly oscillating residual histories. In contrast, H-BACA with $d=32$ achieves significantly better accuracies for most data points particularly as $n_b$ increases. For the EFIE3D (Fig. \ref{fig:accmerge}c) and Frontal3D (Fig. \ref{fig:accmerge}d) matrices, H-BACA with $d=32$ achieves comparable accuracies as H-BACA with $d=1$ for most data points. Note that $d=32$ is significantly better than $d=1$ when the prescribed residual error is large ($\epsilon=10^{-2}$). This agrees with the residual histories in Fig. \ref{fig:history:em3d} and Fig. \ref{fig:history:frontal} as they lie below the singular value spectra when iteration count $k$ is small. 
	
	\subsection{Efficiency}
	This subsection provides six examples to verify the computational complexity estimates in Table 1. H-BACA with leaf-level ACA ($d=1$) and BACA ($d=8,16,32,64,128$) is tested for the following matrices: one Gaussian-SUSY matrix with \ylrev{$n=50000$}, $h=1.0$, $\epsilon=10^{-2}$, one Gaussian-MNIST matrix with $n=5000$, $h=3.0$, $\epsilon=10^{-2}$, one EFIE3D matrix for a unit sphere with $n=26268$, $\epsilon=10^{-6}$ and $20$ points per wavelength, one Frontal3D matrix with $n=1250$, $\epsilon=10^{-6}$ and 10 points per wavelength, one Polynomial matrix with \ylrev{$n=10000$}, $h=0.2$, $\epsilon=10^{-4}$, \ylrev{and one Product-of-random matrix with $n=2500$, $\epsilon=10^{-4}$}. The corresponding $\epsilon$-ranks are 298, \ylrev{137}, 1488, \ylrev{788}, 450 and \ylrev{1000}, respectively. \ylrev{It can be validated that the hierarchical merge operation attains increasing ranks for the Gaussian, EFIE3D and Frontal3D matrices, and relatively constant ranks for the Polynomial, and Product-of-random matrices. All examples use one process except that the Gaussian-SUSY example uses 16 processes.} The CPU times are measured and plotted in Fig. \ref{fig:many}.
	
	Table I predicts that H-BACA exhibits increasing (with a factor of $\sqrt{n_b}$) and constant time when $s_l$ stays constant and increases, respectively. Note that the rank assumption $s_l\approx r$ leading to the $O(\sqrt{n_b})$ computational overhead may not be fully observed for practical values of $n_b$ and $n$. Given one matrix, $s_l$ may stay approximately constant for a limited number of subdivision levels $l$. For example, $s_l$ stay constant for bottom levels of EFIE3D and Frontal3D matrices, and top levels of Polynomial and Product-of-random matrices. This agrees with the observed scalings (w.r.t $n_b$) in Fig. \ref{fig:many:c} - \ref{fig:many:f}. As a reference, the $O(\sqrt{n_b})$ curves are plotted and only small ranges of $n_b$ exhibit the $O(\sqrt{n_b})$ overhead. For the Gaussian matrices, we even observe non-increasing CPU time w.r.t. $n_b$ when $n_b$ is not too big. (see Fig. \ref{fig:many:a} and \ref{fig:many:b}). 
	
	\ylrev{The effects of varying block size $d$ also deserve further discussions. First, larger block size $d$ can significantly improve the robustness of H-BACA for the Gaussian matrices. For example, H-BACA does not achieve desired accuracies due to premature termination for all data points on the $d=1$ curve in Fig. \ref{fig:many:a} and $d=1,8$ curves in Fig. \ref{fig:many:b}. In contrast, H-BACA with larger $d$ attains desired accuracies. Second, larger block size $d$ results in reduced CPU time for the Polynomial  and Frontal3D matrices due to better BLAS performance (see Fig. \ref{fig:many:d} and \ref{fig:many:e}). For the other tested matrices, no significant performance differences have been observed by changing block size $d$. However, for matrices with ranks $s_0\leq d$, larger $d$ and $n_b$ can introduce significant overheads.}

	\subsection{Parallel Performance} 
	Finally, the parallel performance of the H-BACA algorithm is demonstrated via
	strong scaling studies with the EFIE2D, EFIE3D, Product-of-random and Gaussian matrices with process counts $p=8,...,1024$. For the EFIE2D matrices, $n=160000$ and the
	wavenumbers are chosen such that the $\epsilon$-ranks with $\epsilon=10^{-4}$
	are $937$ and $107$, respectively. For the EFIE3D matrices for a unit square, $n=21788$ and the
	wavenumbers are chosen such that the $\epsilon$-ranks with $\epsilon=10^{-6}$
	are $1007$ and $598$, respectively. For the Product-of-random matrices, $n=10000$ and the inner dimension of the product is set to $r=2000$ and $800$, respectively. For the Gaussian matrices with a randomly generated dataset of dimension $50$ and $n=10000$, we choose $h=1.0$ and $h=1.6$ such that the $\epsilon$-ranks with $\epsilon=10^{-3}$ are $2106$ and $191$, respectively. In all examples, the block size and number of leaf-level subblocks in H-BACA are chosen as $d=8$ and $\sqrt{n_b}=\lceil\sqrt{p}\rceil$. The ScaLAPACK \ylrev{block} size is set to $64\times64$. As the reference, we compare to a straightforward parallel implementation of the baseline ACA algorithm which essentially parallelize every operation in ACA with collective MPI communications.  
	
	For all examples, the parallel ACA algorithm stops scaling when $p$ is sufficiently large (see Fig. \ref{fig:ss}). In contrast, the proposed parallel H-BACA algorithm scales up to $p=1024$. In most examples, H-BACA achieves better parallel efficiencies with larger ranks due to better process utilization during the hierarchical merge operation. We also note that ACA outperforms H-BACA for the Product-of-random matrices with small process count $p$ (and $n_b$). This is partially attributed to the $O(\sqrt{n_b})$ overhead observed in Fig. \ref{fig:many:f}.   
	
	Overall, the parallel H-BACA algorithm can achieve reasonably good parallel performances for rank-deficient matrices with modest to large numerical ranks. Not surprisingly, the parallel runtime is dominated by that of ScaLAPACK computation and possible redistributions between each re-compression step as analyzed in Section \nameref{algo_ana}. Also note that the leaf-level BACA compression is embarrassingly parallel for all test cases.

	\section{Conclusion}
	This paper presents a parallel and purely algebraic ACA-type matrix decomposition algorithm given that any matrix entry can be evaluated in $O(1)$ time. Two proposed strategies, BACA and H-BACA, are leveraged to improve the robustness and parallel efficiency of the (baseline) ACA algorithm for general rank-deficient matrices. 
	
	First, the BACA algorithm searches for blocks of row/column pivots via column-pivoted QR on the column/row submatrices at each iteration. The blocking nature of BACA provides a closer estimation of the true residual error and reduces the chance of selecting smaller pivots when compared to ACA. Therefore, BACA exhibits a much smoother and more reliable convergence history. Moreover, blocked operations also benefit from higher flop performance compared to non-blocked ones. For a rank-deficient matrix with dimension $n$ and $\epsilon$-rank $r$, the computational cost of BACA is $O(nr^2)$ assuming the block size constant and iteration count $O(r)$.          
	
	Second, the H-BACA algorithm divides the matrix into $n_b$ similar-sized submatrices each compressed with BACA and then hierarchically merges the results using low-rank arithmetic. Depending on the rank behaviors of submatrices during the merge, the H-BACA may have a computational overhead of $O(\sqrt{n_b})$ yielding the overall computational cost at most $O(nr^2\sqrt{n_b})$. The H-BACA algorithm can be parallelized with distributed-memory machines by assigning each process to one submatrix and leveraging PBLAS and ScaLAPACK for the hierarchical merge operation. Such parallelization strategy yields a much more favorable communication cost when compared to the straightforward parallelization of ACA/BACA with collective MPI routines. Not surprisingly, good parallel performance can be achieved for matrices with modest to large numerical ranks which increases process utilization for each merge operation.   
	
	In contrast to the baseline ACA algorithm, the proposed algorithms exhibit improved robustness and favorable parallel performance with low computational overheads for broad ranges of matrices arising from many science and engineering applications.     
	
	\begin{dci}
		The author(s) declared no potential conflicts of interest with respect to the research, authorship, and/or publication of this article.
	\end{dci}
	
	\begin{funding}
		The author(s) disclosed receipt of the following financial
		support for the research, authorship, and/or publication of
		this article: This research was supported in part by the Exascale Computing Project
		(17-SC-20-SC), a collaborative effort of the U.S. Department of Energy
		Office of Science and the National Nuclear Security Administration,
		and in part by the U.S. Department of Energy, 
		Office of Science, Office of Advanced Scientific Computing Research,
		Scientific Discovery through Advanced Computing (SciDAC) program through the
		FASTMath Institute under Contract No. DE-AC02-05CH11231 at
		Lawrence Berkeley National Laboratory.
	\end{funding}
	
	\begin{acks}
		This research used resources of the National Energy Research Scientific 
		Computing Center (NERSC), a U.S. Department of Energy Office of Science User 
		Facility operated under Contract No. DE-AC02-05CH11231.
	\end{acks}

%	\bibliographystyle{SageH}
%	\bibliography{./baca}

	\begin{biogs}
		\textit{Yang Liu} is a research scientist in the Scalable Solvers Group of the Computational Research Division at Lawrence Berkeley National Laboratory, in Berkeley, California. Dr. Liu received the Ph.D. degree in electrical engineering from the University of Michigan in 2015. From 2015 to 2017, he worked as a postdoctoral fellow at the Radiation Laboratory, University of Michigan. From 2017 to 2019, he worked as a postdoctoral fellow at Lawrence Berkeley National Laboratory, in Berkeley, California. His main research interest is in computational electromagnetics (including fast time-domain integral equation solvers, fast direct integral and differential equation solvers, and multi-physics modeling), numerical linear and multi-linear algebras (including sparse solvers, randomized low-rank, butterfly and tensor algebras), and high-performance scientific computing. Dr. Liu authored and co-authored the Sergei A. Schelkunoff Transactions Prize Paper, APS 2018, second place student paper, ACES 2012, and the first place student paper, FEM 2014.
		
		\noindent\textit{Wissam Sid-Lakhdar} is a postdoctoral researcher in the Scalable Solvers Group of the Computational Research Division at Lawrence Berkeley National Laboratory (LBL). He is currently working on the development of autotuning algorithms and software supported by the Exascale Computing Project (ECP). Before joining LBL, he was a postdoctoral fellow in the PARASOL laboratory in the Computer Science and Engineering Department of Texas A\&M University. He was then working on autotuning batched QR factorization kernels on GPUs. He obtained his Ph.D. at Ecole Normale Superieur of Lyon, where his work targeted the scalability of sparse linear algebra methods on heterogeneous architectures.
		
		\noindent\textit{Elizaveta Rebrova} is an Assistant Adjunct Professor in the Department of Mathematics of the University of California in Los Angeles. In Summers 2017 and 2018 she worked on the interplay between machine learning and numerical linear algebra in the Scalable Solvers Group of the Computational Research Division at Lawrence Berkeley National Laboratory. She received the Ph.D. degree in mathematics from the University of Michigan in 2018. Her main research interests are high-dimensional probability and random matrix theory, and their applications to high-dimensional data science and linear algebra.
			
		\noindent\textit{Pieter Ghysels} is a research scientist in the Scalable Solvers Group of the Computational Research Division at Lawrence Berkeley National Laboratory, in Berkeley, California. His main interests are in High Performance Computing (HPC) and linear algebra. Pieter has expertise in both iterative methods and direct methods for the solution of systems of linear equations. He is the main developer of the STRUMPACK software library which offers a direct solver and preconditioners for large sparse linear systems as well as memory efficient representations of structured dense matrices.  Pieter Ghysels received an engineering degree (in 2006) and completed a PhD in engineering Sciences, both at the (Flemmish) Catholic University in Leuven, Belgium. From 2010-2013, Pieter worked at the Universiteit Antwerpen (University of Antwerp, Belgium) and at the Intel Exascience Lab Flanders.
		
		\noindent\textit{Xiaoye Sherry Li} is a Senior Scientist in the Computational Research Division, Lawrence Berkeley National Laboratory. She has worked on diverse problems in high performance scientific computations, including parallel computing, sparse matrix computations, high precision arithmetic, and combinatorial scientific computing. She has (co)authored over 110 publications, and contributed to several book chapters. She is the lead developer of SuperLU, a widely-used sparse direct solver, and has contributed to the development of several other mathematical libraries, including ARPREC, LAPACK, PDSLin, STRUMPACK, and XBLAS. She has collaborated with many domain scientists to deploy the advanced mathematical software in their application codes, including those from accelerator engineering, chemical science, earth science, plasma fusion energy science, and materials science. She earned Ph.D. in Computer Science from UC Berkeley in 1996. She has served on the editorial boards of the SIAM J. Scientific Comput. and ACM Trans. Math. Software, as well as many program committees of the scientific conferences. She is a SIAM Fellow and an ACM Senior Member.

	\end{biogs}

\end{document}